\documentclass[preprint,11pt]{elsarticle}

\usepackage{amssymb}
\usepackage{amsmath}
\usepackage{amsfonts}
\usepackage{amssymb}
\usepackage{indentfirst,latexsym,bm}
\usepackage{amsthm}
\usepackage{color,xcolor}
\usepackage[all]{xy}
\input{mathrsfs.sty}
\newtheorem{theorem}{Theorem}[section]
\newtheorem{lemma}[theorem]{Lemma}
\newtheorem{corollary}[theorem]{Corollary}
\newtheorem{proposition}[theorem]{Proposition}
\newtheorem{question}[theorem]{Question}

\theoremstyle{definition}
\newtheorem{definition}{Definition}[section]

\theoremstyle{definition}

\theoremstyle{remark}
\newtheorem{remark}{Remark}[section]

\numberwithin{equation}{section}

\journal{XXX}

\begin{document}

\begin{frontmatter}



\title{The generalized polar decomposition, the weak complementarity and the parallel sum for adjointable operators on Hilbert $C^*$-modules}
\author[shnu]{Xiaofeng Zhang}
\ead{xfzhang8103@163.com}
\author[shnu]{Xiaoyi Tian}
\ead{tianxytian@163.com}
\author[shnu]{Qingxiang Xu}
\ead{qingxiang$\_$xu@126.com}
\address[shnu]{Department of Mathematics, Shanghai Normal University, Shanghai 200234, PR China}

\begin{abstract}This paper deals mainly with some aspects of the adjointable operators on Hilbert $C^*$-modules. A new tool called the generalized polar decomposition for each adjointable operator is introduced and clarified. As an application, the general theory of the weakly complementable operators is  set up in the framework of Hilbert $C^*$-modules. It is proved that there exists an operator equation which has a unique solution,  whereas this unique solution fails to be the reduced solution.
Some investigations are also carried out in the Hilbert space case. It is proved that there exist a closed subspace $M$ of certain Hilbert space $K$ and  an operator $T\in \mathbb{B}(K)$ such that $T$ is $(M,M)$-weakly complementable, whereas $T$ fails to be $(M,M)$-complementable. The solvability of the equation
  $$A:B=X^*AX+(I-X)^*B(I-X) \quad (X\in\mathbb{B}(H))$$  is also dealt with in the Hilbert space case,
where $A,B\in \mathbb{B}(H)$ are two general positive operators, and $A:B$ denotes their parallel sum.  Among other things, it is shown that there exist certain positive operators $A$ and $B$ on the Hilbert space $\ell^2(\mathbb{N})\oplus \ell^2(\mathbb{N})$ such that the above equation has no solution.
\end{abstract}

\begin{keyword}Hilbert $C^*$-module, generalized polar decomposition, complementable operator,  reduced solution, parallel sum
\MSC 46L08, 47A05



\end{keyword}

\end{frontmatter}


\section{Introduction}

Hilbert $C^*$-module is one of the natural generalizations of Hilbert space, whose inner-product is taken in a (general) $C^*$-algebra $\mathfrak{A}$ instead of the complex field $\mathbb{C}$. Usually, the underlying $C^*$-algebra $\mathfrak{A}$ is not so simple as $\mathbb{C}$,  therefore many new phenomena may happen when one deals with operators on Hilbert $C^*$-modules. For instance, a closed submodule of a Hilbert $C^*$-module may fail to be orthogonally complemented \cite[P.\,7]{Lance}, so an adjointable operator between two Hilbert $C^*$-modules may have no polar decomposition \cite[Example~3.15]{LLX}, and the Douglas' theorem known for operator equation \cite[Theorem~1]{Douglas} can be not applicable in the Hilbert $C^*$-module case. However in the Hilbert space case, the polar decomposition and the Douglas' theorem usually serve as the main tools and are frequently used in the study of such topics as the parallel sum of two operators \cite{Fillmore-Williams} and its generalizations including the shorted operators \cite{Anderson,AT}, the weak parallel sum and the bilateral shorted operators \cite{ACS}.

We recall briefly some basic knowledge about adjointable operators on  Hilbert $C^*$-modules. The reader is referred to
\cite{Lance,MT,Paschke} for more details; see also \cite{Frank,Frank02,KS,Manuilov,MT02} for some interesting characterizations of Hilbert $C^*$-modules. As is mentioned before, Hilbert $C^*$-modules are generalizations of Hilbert spaces by allowing the inner-product to take values in a $C^*$-algebra rather than in the complex field. Given Hilbert modules $H$ and $K$ over a $C^*$-algebra $\mathfrak{A}$, let $\mathcal{L}(H,K)$ denote the set of all adjointable operators from $H$ to $K$.
It is known that every adjointable operator $A\in\mathcal{L}(H,K)$ is a bounded linear operator, which is also $\mathfrak{A}$-linear.
When $A\in\mathcal{L}(H,K)$, let $\mathcal{R}(A)$, $\mathcal{N}(A)$, $A^*$, $|A|$ and $A|_M$ denote the range of $A$, the null space of $A$, the adjoint operator of $A$, the square root of $A^*A$ and  the restriction of $A$ on a submodule $M$ of $H$, respectively.
  If $H=K$, the notation $\mathcal{L}(H,H)$ is simplified to  $\mathcal{L}(H)$, and we use $\sigma(A)$ to denote the spectrum of $A$ with respect to $\mathcal{L}(H)$.
 Let $I_H$ or simply $I$ denote the identity operator on $H$.  The set of all self-adjoint (resp.\,positive) elements in the unital $C^*$-algebra $\mathcal{L}(H)$ is designated by $\mathcal{L}(H)_{\mbox{sa}}$ (resp.\,$\mathcal{L}(H)_+$). When $A\in\mathcal{L}(H)_+$, the notation $A\ge 0$ is also used to indicate that $A$ is a positive operator on $H$. An operator in $\mathcal{L}(H)$ is referred to be positive definite if it is positive and is invertible. In the special case that $H$ is a Hilbert space, we use the notation $\mathbb{B}(H)$ instead of $\mathcal{L}(H)$.

Suppose that $H$ is a Hilbert module over a $C^*$-algebra $\mathfrak{A}$. An operator $P\in \mathcal{L}(H)$ is called an idempotent if
$P^2=P$. If furthermore $P^*=P$, then $P$ is called a projection. A closed submodule $M$ of  $H$ is said to be orthogonally complemented  in $H$ if $H=M+ M^\bot$, where
$M^\bot=\big\{x\in H:\langle x,y\rangle=0\ \mbox{for every}\ y\in
M\big\}$. In this case, the projection from $H$ onto $M$ is denoted by $P_M$. We use the notations  ``$\oplus$" and ``$+$" with different meanings. For Hilbert
$\mathfrak{A}$-modules $H_1$ and $H_2$, let $H_1\oplus
H_2$ be the Hilbert $\mathfrak{A}$-module defined by $$H_1\oplus
H_2=\left\{\binom{h_1}{h_2}:h_i\in H_i, i=1,2\right\},$$ whose $\mathfrak{A}$-valued inner-product is given by $$\left\langle\binom{h_1}{h_2},\binom{h_1'}{h_2'}\right\rangle=\langle h_1,h_1' \rangle + \langle h_2,h_2' \rangle$$ for $h_i,h_i' \in H_i,i=1,2$. If $M$ and $N$
are non-empty subsets of a Hilbert $\mathfrak{A}$-module, we denote their algebraic sum by $M+N$. That is, $$M+
N=\big\{x+y:x\in M, y\in N \big\}.$$

Unless otherwise specified, throughout the rest of this paper $\mathbb{N}$ is the set of natural numbers, $\mathbb{C}$ is the complex field, $\mathfrak{A}$ is a $C^*$-algebra, $E,H$ and $K$ are Hilbert modules over $\mathfrak{A}$. The main purpose of this paper is, in the framework of the adjointable operators on Hilbert $C^*$-modules, to set up the general theory of the bilateral shorted operators  by following the line initiated in \cite{ACS}. We begin with the introduction of a new tool called the generalized polar decomposition. For every $T\in\mathcal{L}(H,K)$ and every number $\alpha\in (0,1)$, the generalized polar decomposition of $T$ with respect to the parameter $\alpha$ is provided in Theorem~\ref{thm:T=U|T|} by following the techniques employed in the proofs of \cite[Lemma~1.4.4 and Proposition~1.4.5]{Pedersen}. Such a new tool is particularly
useful in dealing with some inequalities related to Hilbert $C^*$-modules. For instance, it can be used to deal with the generalized version of the mixed Schwarz inequality \cite[Theorem~3.2]{SMXZ}. Furthermore, by utilizing this new tool, it will be shown in Corollary~\ref{cor:T and Vt}  that for every $T\in\mathcal{L}(H,K)$,
system \eqref{red sys T} has the  reduced solution denoted by $V_T$, while $V_T^*$ is exactly the reduced solution of another system \eqref{red sys T star}.
This makes it possible to introduce the bilateral shorted operators in the setting of Hilbert $C^*$-modules.

Given Hilbert modules $H$ and $K$ over a $C^*$-algebra $\mathfrak{A}$, let $M$ and $ N$ be orthogonally complemented closed submodules of $H$ and $K$, respectively.
 Up to unitary equivalence, $T\in \mathcal{L}(H,K)$ can be partitioned as \eqref{2 form A} with the operator $T_{22}$ in its lower right corner. Based on the induced operators $V_{T_{22}}$ and $V_{T_{22}}^*$ together with $|T_{22}|^\frac12$ and $|T_{22}^*|^\frac12$, conditions for $T$ to be $(M,N)$-weakly complementable are figured out in Definition~\ref{defn:(S,T)-weakly}, along with a new formula \eqref{formula with half} for the bilateral shorted operator $T_{/(M,N)}$ in the case that  $T$ is  $(M,N)$-weakly complementable. As shown by Propositions~\ref{prop:polar decomposition-weak reduced} and \ref{prop:s to w-01}, this operator $T_{/(M,N)}$ is the same as that considered in \cite{ACS} whenever $T$  has the polar decomposition or $T$ is
$(M,N)$-complementable. Thus, the term of the bilateral shorted operator has been generalized from the Hilbert space case to the Hilbert $C^*$-module case.
Although the contents of Propositions~\ref{cor:range A(M,N)} and Corollary~\ref{cor:range characterization00} are almost the same as that of \cite[Corollaries 4.6 and 4.9]{ACS},
the proofs are however somehow different. For, the proofs of
\cite[Corollaries 4.6 and 4.9]{ACS} are based on
\cite[Lemma~4.4 and Proposition~4.5]{ACS}, which are figured out by using some well-known knowledge on a Hilbert space $H$.  For instance, the weak compactness of the closed unital ball of $H$ (a fact derived from  Riesz theorem together with Alaoglu-Bourbaki theorem) and the polar decomposition for each operator in $\mathbb{B}(H)$  are served as the tools in the proof of \cite[Lemma~4.4]{ACS}.
It is notable that these tools can not always be applicable in the Hilbert $C^*$-module case.

 Let $A\in\mathcal{L}(H, K)$ and $C\in\mathcal{L}(E,K)$. In the special case that all of $E,H$ and $K$ are Hilbert spaces, it can be inferred directly from the Douglas' theorem that system \eqref{equ:system A X equals A prime} is solvable iff $\mathcal{R}(C)\subseteq \mathcal{R}(A)$, which in turn, holds iff system \eqref{equ:system A X equals A prime} has a reduced solution. When we go back to the Hilbert $C^*$-module case, it is known that the inclusion $\mathcal{R}(C)\subseteq \mathcal{R}(A)$ generally  can not ensure the solvability  of \eqref{equ:system A X equals A prime}.
Actually, this makes the main difference between our characterization for $(M,N)$-complementable operators (Proposition~\ref{prop of S T complementable}) and that carried out in \cite[Section~3]{ACS}. To the best of our knowledge,  the answer to the following question seems to be unknown in the setting of Hilbert $C^*$-module operators:

\begin{question}Is it true that every solvable operator equation has the reduced solution?
\end{question}
In this paper, we will give the negative answer to the above question; see Theorem~\ref{thm:solvability does not imp reduced solu} for the details.

Although the parallel sum for matrices and operators  has been intensely \-studied \cite{Anderson-Duffin,ACS,FMXZ,Hansen,LSX,Mitra-Odell,Morley,TWD}, there are still some issues remain to be unknown. Suppose that $H$ is a Hilbert space and  $A,B\in \mathbb{B}(H)$ are both  positive  definite. It can be inferred from \cite[Lemma~18]{Anderson-Duffin} that
$$\big\langle (A:B)x, x\big\rangle\le \big\langle Ay,y\big\rangle+\big\langle B(x-y),x-y\big\rangle$$
for every $x,y\in H$. Evidently, the inequality above can be rephrased as
 \begin{equation}\label{equ: A:B geq C^*AC+(I-C)*B(I-C)}
  A:B\le C^*AC+(I-C)^*B(I-C), \quad \forall\, C \in \mathbb{B}(H).
\end{equation}
 For a proof of the latter inequality, the reader is referred to \cite[Theorem~2.1]{Hansen} (see also Lemma~\ref{thm:1st on parallel sum} in this paper).
 Replacing $A$ and $B$ with $A+\frac1n I$ and $B+\frac1n I$ respectively, we see that \eqref{equ: A:B geq C^*AC+(I-C)*B(I-C)}
 is also true for each pair of positive operators $A$ and $B$ in $\mathbb{B}(H)$. As far as we know, little has been done in the literature  on the  solvability  of the following equation
 \begin{equation}\label{equality is obtained}
  A:B=X^*AX+(I-X)^*B(I-X), \quad X\in\mathbb{B}(H),
\end{equation}
where $A,B\in \mathbb{B}(H)$ are two general positive operators. An investigation of the  solvability  of the above equation is undertaken in this paper. Among other things, it is shown in Theorem~\ref{thm:unsolvable} that  equation \eqref{equality is obtained} can be unsolvable for some positive operators $A$ and $B$ on the Hilbert space $\ell^2(\mathbb{N})\oplus \ell^2(\mathbb{N})$.

 The rest of the paper is organized as follows. Some basic knowledge about the adjointable operators on Hilbert $C^*$-modules are provided in Section~\ref{sec:pre}. Section~\ref{sec:generalized polar decomposition} is devoted to the construction of the generalized polar decomposition for the adjointable operators, while Sections~\ref{sec:complementarity} and \ref{sec:weak complementarity} are focused on the study of the $(M,N)$-complementable operators and the $(M,N)$-weakly complementable operators, respectively. Two special positive operators on the Hilbert space $\ell^2(\mathbb{N})\oplus \ell^2(\mathbb{N})$ are constructed and are dealt with in Section~\ref{subsec:two positive operators}. As applications, some investigations around these two special positive operators are carried out; see Sections~\ref{subsec:complementarity}--\ref{subsec:equation wrt parallel sum} for the details.
\section{Preliminaries}\label{sec:pre}

Given a projection $P\in\mathcal{L}(H)$ such that $P\ne 0$ and $P\ne I$,  let $U_P: H\to \mathcal{R}(P)\oplus \mathcal{N}(P)$ be the unitary defined  by
\begin{equation}\label{equ:unitary operator induced by P}U_Ph=\binom{Ph}{(I-P)h}\quad (h\in H).
\end{equation}
It is clear that for every $h_1\in\mathcal{R}(P)$ and $h_2\in\mathcal{N}(P)$,
\begin{equation*}\label{expression of inverse of U P}U_P^{*}\binom{h_1}{h_2}=h_1+h_2.\end{equation*}
For every $T\in \mathcal{L}(H)$,  we have  (see e.g.\,\cite[Section~2.2]{Xu-Wei-Gu})
 \begin{equation*}
U_PTU_P^{*}=\left(
            \begin{array}{cc}
              PTP|_{\mathcal{R}(P)} & PT(I-P)|_{\mathcal{N}(P)} \\
              (I-P)TP|_{\mathcal{R}(P)} &(I-P)T(I-P)|_{\mathcal{N}(P)}\\
            \end{array}
          \right).
\end{equation*}
Specifically,
 \begin{equation*}\label{equ:block matrix P}
U_PPU_P^{*}=\left(
            \begin{array}{cc}
              I_{\mathcal{R}(P)} & 0 \\
              0 &0\\
            \end{array}
          \right),\quad U_P(I-P)U_P^{*}=\left(
            \begin{array}{cc}
              0 & 0 \\
              0 &I_{\mathcal{N}(P)}\\
            \end{array}
          \right).
\end{equation*}
Conversely, given every $X=\left(\begin{array}{ccc} X_{11}&
X_{12}\\X_{21}&X_{22}\end{array}\right)\in \mathcal{L}\big(\mathcal{R}(P)\oplus \mathcal{N}(P)\big)$ and every $h\in H$, we have
\begin{align*}U_P^{*}XU_P(h)=&U_P^{*}X\binom{Ph}{h-Ph}=(X_{11}+X_{21})Ph+(X_{12}+X_{22})(I-P)h.
\end{align*}
It follows that
\begin{equation*}U_P^{*}XU_P=(X_{11}+X_{21})P+(X_{12}+X_{22})(I-P).\end{equation*}

\begin{remark}
It is notable that when $P=0$ or $P=I$, all the derivations above are also valid.
\end{remark}

\begin{definition}\label{defn of the reduced solution}{\rm \cite[Section~1]{MMX}}
Consider the system
\begin{equation}\label{equ:system A X equals A prime}  AX=C, \quad X\in\mathcal{L}(E,H),\end{equation}
where $A\in\mathcal{L}(H, K)$ and $C\in\mathcal{L}(E,K)$ are given. An operator $D\in\mathcal{L}(E,H)$ is said to be the reduced solution of \eqref{equ:system A X equals A prime} if $$AD=C\quad\mbox{and}\quad \mathcal{R}(D)\subseteq \overline{\mathcal{R}(A^*)}.$$
\end{definition}
Observe that $\overline{\mathcal{R}(A^*)}\subseteq \mathcal{N}(A)^\bot$, so if such a reduced solution exists, then it is unique.

\begin{definition}\label{defn:defn of polar decomposition}{\rm\cite[Definition~3.10]{LLX}} For each $T$ in $\mathcal{L}(H,K)$, its polar decomposition (if it exists) is represented by
\begin{equation}\label{equ:two conditions of polar decomposition}T=U|T|\quad \mbox{and}\quad U^*U=P_{\overline{\mathcal{R}(T^*)}},\end{equation}
where $U\in\mathcal{L}(H,K)$ is a partial isometry.
\end{definition}

It is known (see e.g.\,\cite[Lemma~3.9]{LLX}) that there exists at most one partial isometry $U\in\mathcal{L}(H,K)$ satisfying  \eqref{equ:two conditions of polar decomposition}.

Next, we provide a couple of lemmas, which will be used in the sequel.
\begin{lemma}\label{lem:Range closure of TT and T and closures are the same-01} {\rm \cite[Proposition~3.7]{Lance}} $\overline {\mathcal{R}(TT^*)}=\overline{ \mathcal{R}(T)}$ for every  $T\in\mathcal{L}(H,K)$.
\end{lemma}

\begin{lemma}\label{lem:Range closure of TT and T and closures are the same-02}
{\rm (\cite[Proposition~2.9]{LLX} and \cite[Lemma~2.2]{VMX})}
$\overline{\mathcal{R}(T^\alpha)}=\overline{\mathcal{R}(T)}$ and $\mathcal{N}(T)=\mathcal{N}(T^{\alpha})$ for every  $T\in \mathcal{L}(H)_+$ and $\alpha>0$.
\end{lemma}

\begin{lemma}\label{lem:dougla's theorem}{\rm (cf.\,\cite[Theorem~3.2]{Fang-Moslehian-Xu})} For every $A\in\mathcal{L}(H, K)$, the following statements are equivalent:
\begin{enumerate}
\item[{\rm (i)}] $\overline{\mathcal{R}(A^*)}$ is orthogonally complemented in $H$;
\item[{\rm (ii)}]Given every $C\in\mathcal{L}(E,K)$ with $\mathcal{R}(C)\subseteq \mathcal{R}(A)$, system  \eqref{equ:system A X equals A prime}  has the reduced solution.
\end{enumerate}
\end{lemma}

\begin{lemma}{\rm \cite[Lemma~3.6 and Theorem~3.8]{LLX}}\label{lem:polar decomposition of T star}
For every $T\in\mathcal{L}(H,K)$, the following statements are equivalent:
\begin{enumerate}
\item[{\rm (i)}] $T$ has the polar decomposition;
\item[{\rm (ii)}]$T^*$ has the polar decomposition;
\item[{\rm (iii)}]$\overline{\mathcal{R}(T)}$ and $\overline{\mathcal{R}(T^*)}$ are orthogonally complemented in $K$ and $H$, respectively.
\end{enumerate}
If {\rm (i)--(iii)} are satisfied and $T$ has the polar decomposition represented by \eqref{equ:two conditions of polar decomposition}, then
the polar decomposition of $T^*$ can be expressed as
\begin{equation}\label{equ:the polar decomposition of T star-pre stage}T^*=U^*|T^*|\quad \mbox{and}\quad UU^*=P_{\overline{\mathcal{R}(T)}}.
\end{equation}
 \end{lemma}

\begin{remark}
Suppose that $T\in\mathcal{L}(H,K)$ has the polar decomposition represented by \eqref{equ:two conditions of polar decomposition}. For simplicity, henceforth  we just say that $T$ has the polar decomposition $T=U|T|$.
\end{remark}

\section{The generalized polar decompositions for adjointable operators}\label{sec:generalized polar decomposition}

As mentioned before, an adjointable operator acting on a Hilbert $C^*$-module may fail to have the polar decomposition. This leads us to introduce the generalized polar decompositions  as follows.
\begin{theorem}\label{thm:T=U|T|} For every $T\in\mathcal{L}(H,K)$ and $\alpha\in (0, 1)$, there exists $U\in\mathcal{L}(H,K)$ such that
\begin{enumerate}
\item[{\rm (i)}] $T=U|T|^\alpha$ and $T^*=U^*|T^*|^\alpha$;
\item[{\rm (ii)}]$U^*U=|T|^{2(1-\alpha)}$ and $UU^*=|T^*|^{2(1-\alpha)}$;
\item[{\rm (iii)}] $U|T|^\beta=|T^*|^\beta U$ for every $\beta>0$.
\end{enumerate}
\end{theorem}
\begin{proof} Let $\alpha\in (0,1)$ be chosen arbitrarily. First, we consider the special case that $H=K$. For each $n\in\mathbb{N}$, let
$$U_n= T\left[1/n\cdot I+T^*T\right]^{-\frac{1}{2}}(T^*T)^\frac{1-\alpha}{2}.$$
 By the proof of \cite[Proposition~1.4.5]{Pedersen},  $\{U_n\}_{n=1}^\infty$ converges to some $U\in\mathcal{L}(H)$ in the norm topology such that
  $$T=U(T^*T)^{\frac{\alpha}{2}}=U|T|^\alpha.$$
It is clear that
\begin{equation}\label{equ:u n star u n}U_n^*U_n=(T^*T)^{2-\alpha}\left[1/n\cdot I+T^*T\right]^{-1}.\end{equation}
Let $M=\|T\|^2$ and define
$$f_n(t)=\frac{t^{2-\alpha}}{1/n+t},\quad f(t)=t^{1-\alpha}\quad (t\in [0, M]).$$
Obviously, $\{f_n\}_{n=1}^\infty$  converges to $f$ pointwisely, $f_n (n\in\mathbb{N})$ and $f$ are all continuous on $[0, M]$ such that for each $t\in [0, M]$, the sequence $\{f_n(t)\}_{n=1}^\infty$ is monotonically increasing. It follows from Dini's theorem that $\{f_n\}_{n=1}^\infty$ converges uniformly to $f$ on $[0, M]$. Therefore, it can be concluded from \eqref{equ:u n star u n}
 that
 $$U^*U=(T^*T)^{1-\alpha}=|T|^{2(1-\alpha)}.$$

A simple replacement of $T$ with $T^*$ yields $T^*=V|T^*|^\alpha$ and
$V^*V=|T^*|^{2(1-\alpha)}$ such that $V=\lim\limits_{n\to\infty}V_n$, where
$$V_n= T^*\left[1/n\cdot I+TT^*\right]^{-\frac{1}{2}}(TT^*)^\frac{1-\alpha}{2}.$$
Observe that for each $n\in\mathbb{N}$,
\begin{align*}V_n=&\left[1/n\cdot I+T^*T\right]^{-\frac{1}{2}}T^*(TT^*)^\frac{1-\alpha}{2}\\
=&\left[1/n\cdot I+T^*T\right]^{-\frac{1}{2}}(T^*T)^\frac{1-\alpha}{2}T^*=U_n^*.
\end{align*}
Taking the limits as $n\to\infty$ yields $V=U^*$. This shows the validity of (i) and (ii).

For every $\beta>0$, we have
$T^*|T^*|^\beta=|T|^\beta T^*,$ so
\begin{align*}
U^*|T^*|^{\beta}=&\lim\limits_{n\to \infty} U_n^*|T^*|^{\beta}=\lim\limits_{n\to \infty}\left[1/n\cdot I+T^*T\right]^{-\frac{1}{2}}(T^*T)^\frac{1-\alpha}{2}T^*|T^*|^{\beta}\\
=&\lim\limits_{n\to \infty}\left[T\left[1/n\cdot I+T^*T\right]^{-\frac{1}{2}}(T^*T)^{\frac{1-\alpha}{2}}|T|^{\beta}\right]^*\\
 =&\lim_{n\to\infty}(U_n|T|^\beta)^*=(U|T|^\beta)^*=|T|^\beta U^*.
\end{align*}
Hence, $|T^*|^\beta U=U|T|^\beta$ by taking $*$-operation.

Next, we consider the general case that $H\ne K$. Let $\widetilde{T}\in\mathcal{L}(H\oplus K)$ be defined by
\begin{equation*}\widetilde{T}=\left(
                                \begin{array}{cc}
                                  0 & 0 \\
                                  T & 0 \\
                                \end{array}
                              \right).
\end{equation*}
It is clear that
$$|\widetilde{T}|=\left(
                    \begin{array}{cc}
                      |T| & 0 \\
                      0 & 0 \\
                    \end{array}
                  \right),\quad |\widetilde{T}^*|=\left(
                                                    \begin{array}{cc}
                                                      0 & 0 \\
                                                      0 & |T^*| \\
                                                    \end{array}
                                                  \right).$$
By the previous argument, there exists
$$\widetilde{U}=\left(
                  \begin{array}{cc}
                    U_{11} & U_{12} \\
                    U & U_{22} \\
                  \end{array}
                \right)\in\mathcal{L}(H\oplus K)$$ with $U\in\mathcal{L}(H,K)$ such that items (i)--(iii) are satisfied with $T$ and $U$ therein be replaced by $\widetilde{T}$ and $\tilde{U}$, respectively. It is routine to verify that $\tilde{U}=\begin{pmatrix}0&0\\
                U&0\end{pmatrix}$ and items (i)--(iii) are all satisfied.
\end{proof}

As a consequence of the preceding theorem, we provide a technical result as follows.
\begin{corollary}\label{cor:T and Vt}For every $T\in\mathcal{L}(H,K)$, there exists uniquely an operator $V_T\in \mathcal{L}(H,K)$ such that $V_T$ and $V_T^*$ are respectively the reduced solutions of the systems
\begin{align}\label{red sys T}&|T^*|^\frac12 X=T, \quad X\in\mathcal{L}(H,K),\\
\label{red sys T star}&|T|^{\frac{1}{2}}Y=T^*,\quad Y\in\mathcal{L}(K,H).
\end{align}
Moreover, it holds that
\begin{align}&\label{V T T star}V_T^*V_T=|T|,\quad V_TV_T^*=|T^*|,\\
\label{V T commutative}&V_T|T|^\beta=|T^*|^\beta V_T,\quad\forall \beta>0,\\
&\label{range closure of V T}\overline{\mathcal{R}(V_T)}=\overline{\mathcal{R}(T)},\quad \overline{\mathcal{R}(V_T^*)}=\overline{\mathcal{R}(T^*)}.
\end{align}
\end{corollary}
\begin{proof} We may choose any $\alpha \in (\frac12, 1)$ in Theorem~\ref{thm:T=U|T|} to get an operator $U_{\alpha}\in\mathcal{L}(H,K)$ satisfying
\begin{align}\label{prop-01 of A 22}&T=U_{\alpha}|T|^\alpha, \quad T^*=U_{\alpha}^*|T^*|^\alpha,\\
\label{prop-02 of A 22}&U_{\alpha}^*U_{\alpha}=|T|^{2(1-\alpha)},\quad U_{\alpha}U_{\alpha}^*=|T^*|^{2(1-\alpha)},\\
\label{prop-03 of A 22}&U_{\alpha}|T|^\beta=|T^*|^\beta U_{\alpha},\quad \forall \beta>0.
\end{align}
Let $V_T\in \mathcal{L}(H,K)$ be defined by
\begin{equation}\label{equ:auxilirary operator W}
V_T=U_{\alpha}|T|^{\alpha-\frac12}.
\end{equation}
It follows from \eqref{equ:auxilirary operator W}, \eqref{prop-01 of A 22} and \eqref{prop-03 of A 22} that
\begin{align*}&|T|^\frac12 V_T^*=|T|^\alpha U_{\alpha}^*=T^*,\\
&|T^*|^\frac12 V_T=(|T^*|^\frac12 U_{\alpha})|T|^{\alpha-\frac{1}{2}}=(U_{\alpha}|T|^\frac12)|T|^{\alpha-\frac{1}{2}}=T.
\end{align*}
Hence,
\begin{equation}\label{red sys T++}|T^*|^\frac12 V_T=T, \quad |T|^\frac12 V_T^*=T^*.
\end{equation}
In addition, by \eqref{equ:auxilirary operator W}, \eqref{prop-02 of A 22} and \eqref{prop-03 of A 22} we have
\begin{align*}&V_T^*V_T=|T|^{\alpha-\frac12}\cdot U_{\alpha}^*U_{\alpha}\cdot |T|^{\alpha-\frac12}=|T|,\\
&V_TV_T^*=U_{\alpha}\cdot |T|^{2\alpha-1}\cdot U_{\alpha}^*=U_{\alpha}\cdot (U_{\alpha}|T|^{2\alpha-1})^*=U_{\alpha}U_{\alpha}^*\cdot |T^*|^{2\alpha-1}=|T^*|.
\end{align*}
In view of  Lemmas~\ref{lem:Range closure of TT and T and closures are the same-01} and \ref{lem:Range closure of TT and T and closures are the same-02},
we see that \eqref{range closure of V T} is immediate from \eqref{V T T star}.
As a result, we may use \eqref{red sys T++}  and \eqref{range closure of V T} together with Lemmas~\ref{lem:Range closure of TT and T and closures are the same-01} and \ref{lem:Range closure of TT and T and closures are the same-02} to conclude that $V_T$ and $V_T^*$ are the reduced solutions of
\eqref{red sys T} and \eqref{red sys T star}, respectively. Since the reduced solution of \eqref{red sys T} is unique,
the uniqueness of $V_T$ follows.

Finally, for every $\beta>0$ we have
\begin{align*}&V_T|T|^\beta=U_{\alpha}|T|^{{\alpha-\frac12}+\beta}=|T^*|^{{\alpha-\frac12}+\beta}U_{\alpha}=|T^*|^\beta\big(|T^*|^{{\alpha-\frac12}}U_{\alpha}\big)=|T^*|^\beta V_T.
\end{align*}

\end{proof}
\begin{remark}\label{condition:A_22 has polar } It is helpful to derive a formula for $V_T$
when $T\in\mathcal{L}(H,K)$ has the polar decomposition.
Assume that \eqref{equ:two conditions of polar decomposition} is satisfied. Let
\begin{equation}\label{V and T-01}V_T=|T^*|^\frac12 U.\end{equation}
A simple use of \eqref{V and T-01} and \eqref{equ:the polar decomposition of T star-pre stage} gives
$T=|T^*|^\frac12 V_T$. It is obvious that $\mathcal{R}(V_T)\subseteq \mathcal{R}\big(|T^*|^\frac12\big)$,  so $V_T$ is exactly the reduced solution of \eqref{red sys T}.
Furthermore, we have
\begin{equation}\label{V and T-02}V_TU^*=|T^*|^\frac12 UU^*=|T^*|^\frac12.\end{equation}
So, we may combine \eqref{V and T-01} and \eqref{V and T-02} to get $\mathcal{R}(V_T)=\mathcal{R}(|T^*|^\frac12)$.
Replacing $T$ with $T^*$ yields
\begin{equation}\label{V and T-03}V_T^*=|T|^\frac12 U^*,\quad V_T^*U=|T|^\frac12,\quad  \mathcal{R}(V_T^*)=\mathcal{R}(|T|^\frac12).\end{equation}
\end{remark}

\section{$(M,N)$-complementable operators}\label{sec:complementarity}
Following the line in \cite[Section~3]{ACS},  we study $(M,N)$-complemen\-table operators in the setting of Hilbert $C^*$-module operators.
Suppose that $M$ and $ N$ are orthogonally complemented closed submodules of $H$ and $K$, respectively.
Let $U_{P_{M}}: H\to M\oplus M^{\bot}$ and  $U_{P_{N}}: K\to N\oplus N^{\bot}$
be defined by \eqref{equ:unitary operator induced by P} with $P$ therein be replaced with $P_M$ and $P_N$, respectively. For each $T\in \mathcal{L}(H,K)$, it is easily seen that
\begin{equation}\label{2 form A}U_{P_{N}}TU_{P_{M}}^*=\left(
  \begin{array}{cc}
    T_{11} & T_{12} \\
    T_{21} & T_{22} \\
  \end{array}
\right),\end{equation}
where $P_{M^\perp}=I_H-P_{M}$, $P_{N^\perp}=I_K-P_{N}$ and $T_{11}\in\mathcal{L}(M,N), T_{12}\in\mathcal{L}(M^{\perp},N)$, $T_{21}\in\mathcal{L}(M,N^\perp)$ and $T_{22}\in\mathcal{L}(M^\perp,N^\perp)$ are formulated by
\begin{align}\label{equ:A matrix-01}&T_{11}=P_{N}TP_{M}|_M,\quad T_{12}=P_{N}TP_{M^\perp}|_{M^\perp},  \\
   \label{equ:A matrix-02}& T_{21}=P_{N^\perp}TP_{M}|_M,\quad T_{22}=P_{N^\perp}TP_{M^\perp}|_{M^{\perp}}.
\end{align}

\begin{definition}{\rm \cite[Definition~3.1]{ACS}}\label{defn:(S,T)}
Let $M$ and $N$ be orthogonally complemented closed submodules of $H $ and $K$, respectively. An operator $T\in \mathcal{L}(H,K)$ is said to be $(M,N)$-complementable if
there exist $M_r\in\mathcal{L}(H)$ and $M_l\in\mathcal{L}(K)$ such that
  \begin{align}
    &\label{equ:M_s and M_T-01}(I_H-P_M)M_r=M_r, \quad M_l(I_K-P_N)=M_l, \\
    &\label{equ:M_s and M_T-02}(I_K-P_N)TM_r=(I_K-P_N)T, \quad  M_lT(I_H-P_M)=T(I_H-P_M).
 \end{align}
\end{definition}

\begin{proposition}\label{prop of S T complementable} {\rm (cf.\,\cite[Proposition~3.2]{ACS})} Let $M$ and $N$ be orthogonally complemented closed submodules of $H $ and $K$, respectively. For every $T\in \mathcal{L}(H,K)$, the following statements are equivalent:
\begin{enumerate}
  \item[{\rm (i)}]$T$ is $(M,N)$-complementable.

  \item[{\rm (ii)}]Each of the systems
\begin{align}\label{defn:(S,T)complementable-01}
  T_{22}X&=T_{21}, \quad X\in\mathcal{L}(M,M^\perp),\\
\label{defn:(S,T)complementable-02}  T_{22}^*Y&=T_{12}^*, \quad Y\in\mathcal{L}(N,N^\perp),
\end{align}
has a solution, where $T_{ij}(i,j = 1,2)$ are defined by \eqref{equ:A matrix-01} and \eqref{equ:A matrix-02}.
  \item[{\rm (iii)}]There exist idempotents $P\in\mathcal{L}(H)$ and $Q\in\mathcal{L}(K)$ such that
      \begin{equation}\label{4 conditions of idempotents}\mathcal{R}(P^*)=M,\quad \mathcal{R}(TP)\subseteq N,\quad \mathcal{R}(Q)=N,\quad  \mathcal{R}\big((QT)^*\big)\subseteq M.\end{equation}
  \item[{\rm (iv)}]$H$ and $K$ can be decomposed algebraically  as
  \begin{align*}&H=M^\perp+\mathcal{N}\big[(I_K-P_N)T(I_M-C)P_M\big],\\ &K=N^\perp+\mathcal{N}\big[(I_H-P_M)T^*(I_N-D)P_N\big]\end{align*}
  for some $C\in \mathcal{L}(M,M^\perp)$ and $D\in \mathcal{L}(N,N^\perp)$.
  \end{enumerate}
  \end{proposition}
\begin{proof}
(i)$\Longrightarrow$(ii). It is clear that
\begin{equation}\label{matrix form of S T}U_{P_M}P_MU^*_{P_M}=\left(
    \begin{array}{cc}
      I_{M} & 0 \\
      0 & 0 \\
    \end{array}
  \right), \quad
  U_{P_N}P_NU^*_{P_N}=\left(
    \begin{array}{cc}
      I_{N} & 0 \\
      0 & 0 \\
    \end{array}
  \right),
\end{equation}
which lead by \eqref{equ:M_s and M_T-01} to
\begin{equation}\label{matrix form of M S T}U_{P_M}M_rU^*_{P_M}=\left(
                  \begin{array}{cc}
                    0 & 0 \\
                    C & D \\
                  \end{array}
                \right),
\quad
U_{P_N}M_lU^*_{P_N}=\left(
                  \begin{array}{cc}
                    0 & E \\
                    0 & F \\
                  \end{array}
                \right)
\end{equation}
for some $C\in \mathcal{L}(M, M^\perp)$, $D\in \mathcal{L}(M^\perp)$, $E\in \mathcal{L}(N^\perp, N)$ and $F\in \mathcal{L}(N^\perp)$. Rewrite the first equation in \eqref{equ:M_s and M_T-02} as
$$U_{P_N}(I_K-P_N)U_{P_N}^*\cdot U_{P_N}TU_{P_M}^*\cdot U_{P_M} M_rU_{P_M}^*=U_{P_N}(I_K-P_N)U_{P_N}^*\cdot U_{P_N}TU_{P_M}^*.
$$
Substituting \eqref{matrix form of S T}, \eqref{2 form A} and \eqref{matrix form of M S T} into the above equation yields
\begin{equation*}\label{equ:A21=A22C}
  T_{22}C=T_{21} \quad\mbox{and}\quad T_{22}D=T_{22}.
\end{equation*}
Therefore, the operator $C$ is a solution of \eqref{defn:(S,T)complementable-01}.
Similarly, from the second equation in \eqref{equ:M_s and M_T-02} we can obtain
\begin{equation*}
  ET_{22}=T_{12} \quad\mbox{and}\quad FT_{22}=T_{22},
\end{equation*}
which shows that $E^*$ is a solution of \eqref{defn:(S,T)complementable-02}.

(ii)$\Longrightarrow$(iii). Suppose that $C$ and $E$ are (arbitrary) solutions of \eqref{defn:(S,T)complementable-01} and \eqref{defn:(S,T)complementable-02}, respectively.  Let $P\in \mathcal{L}(H)$ and $Q\in \mathcal{L}(K)$ be the idempotents defined by
$$P=U^*_{P_M}\left(
    \begin{array}{cc}
      I_{M} & 0 \\
      -C & 0 \\
    \end{array}
  \right)U_{P_M}
\quad \mbox{and} \quad
Q=U^*_{P_N}\left(
    \begin{array}{cc}
      I_{N} & -E^* \\
      0 & 0 \\
    \end{array}
  \right)U_{P_N}.$$
Evidently, $\mathcal{R}(U_{P_M}P^*U^*_{P_M})=M\oplus\{0\},$ so
$$\mathcal{R}(P^*)=\mathcal{R}(U_{P_M}^*\cdot U_{P_M}P^*U^*_{P_M})=M.$$
A similar reasoning gives $\mathcal{R}(Q)=N$. Moreover, a simple computation yields
\begin{align*}&U_{P_N}TPU^*_{P_M}=\left(
                                           \begin{array}{cc}
                                             T_{11}-T_{12}C &\quad 0 \\
                                             T_{21}-T_{22}C &\quad 0 \\
                                           \end{array}
                                         \right)
                                         =\left(
                                           \begin{array}{cc}
                                             T_{11}-T_{12}C &\quad 0 \\
                                            0 &\quad 0 \\
                                           \end{array}
                                         \right),\\
&U_{P_M}(QT)^*U^*_{P_N}=\left(
                                           \begin{array}{cc}
                                             T_{11}^*-T_{21}^*E &\quad 0 \\
                                             T_{12}^*-T_{22}^*E &\quad 0 \\
                                           \end{array}
                                         \right)
                                         =\left(
                                           \begin{array}{cc}
                                             T_{11}^*-T_{21}^*E &\quad 0 \\
                                            0 &\quad 0 \\
                                           \end{array}
                                         \right),
\end{align*}
 which imply that $\mathcal{R}(TP)\subseteq N$ and $\mathcal{R}\big((QT)^*\big)\subseteq M$.

(iii)$\Longrightarrow$(iv). Suppose that there exist idempotents $P\in\mathcal{L}(H)$ and $Q\in\mathcal{L}(K)$ such that \eqref{4 conditions of idempotents} is satisfied.
     Since $\mathcal{N}(P)=\mathcal{R}(P^*)^\perp=M^\perp$, we have
\begin{equation}\label{P 2*2}
  U_{P_{M}}PU_{P_{M}}^*=\left(
                          \begin{array}{cc}
                            P_{11} & 0 \\
                            -C & 0 \\
                          \end{array}
                        \right)
\end{equation}
for some $P_{11}\in \mathcal{L}(M)$ and $C\in \mathcal{L}(M,M^\perp)$, which leads by $P^2=P$ to
\begin{equation*}
 P_{11}^2=P_{11}\quad\mbox{and}\quad  CP_{11}=C.
\end{equation*}
Combining the above equation with \eqref{P 2*2} yields
$$U_{P_{M}}P(I_M-P_{11})x=U_{P_{M}}PU_{P_{M}}^*\cdot \begin{pmatrix}(I_M-P_{11})x\\0\end{pmatrix}=\left(
                                                                                                    \begin{array}{c}
                                                                                                      0 \\
                                                                                                      0\\
                                                                                                    \end{array}
                                                                                                  \right)
$$ for every $x\in M$, which means that $$(I_M-P_{11})x\in M\cap \mathcal{N}(P)=M\cap M^\perp=\{0\}.$$
Consequently, $P_{11}=I_M$. As a result,
for each $\xi\in H$ we have
\begin{align*}P\xi=&U_{P_M}^*\left(
                              \begin{array}{cc}
                                I_M & 0 \\
                                -C & 0 \\
                              \end{array}
                            \right)U_{P_M}\xi=U_{P_M}^*
                              \begin{pmatrix}
                                P_M\xi\\ -CP_M\xi
                              \end{pmatrix}=(I_M-C)P_M\xi.
\end{align*}
This shows that $P=(I_M-C)P_M$.

In view of $\mathcal{R}(TP)\subseteq N$ and $P^2=P$, we have
$$(I_K-P_N)TP\cdot P=0,$$
which yields
$$\mathcal{R}(P)\subseteq \mathcal{N}\big[(I_K-P_N)TP\big]=\mathcal{N}\big[(I_K-P_N)T(I_M-C)P_M\big].$$
Observe that
$$H=\mathcal{N}(P)+\mathcal{R}(P)=\mathcal{R}(P^*)^\perp +\mathcal{R}(P)=M^\perp +\mathcal{R}(P).$$
So, $H$ can be decomposed as desired. Similarly, the desired decomposition of $K$  can be derived from the last two conditions in \eqref{4 conditions of idempotents}.

(iv)$\Longrightarrow$(i). Let $W\in \mathcal{L}(M,K)$ be defined by
$$W=(I_K-P_N)T(I_M-C).$$
For each $x\in M$, from \eqref{2 form A} it is easy to show that
\begin{align}U_{P_{N}}Wx=& U_{P_{N}}(I_K-P_N)U_{P_{N}}^*\cdot U_{P_{N}} TU_{P_{M}}^*\cdot  U_{P_{M}}(I_M-C)x\nonumber\\
\label{useful-01}=&\left(
                     \begin{array}{cc}
                       0 & 0 \\
                       T_{21} & T_{22} \\
                     \end{array}
                   \right)\left(
                            \begin{array}{c}
                              x \\
                              -Cx \\
                            \end{array}
                          \right)=
\left(
       \begin{array}{c}
         0 \\
         T_{21}x-T_{22}Cx \\
       \end{array}
     \right).\end{align}
Utilizing the  decomposition of $H$ as asserted, we have
$x=x_1+x_2$
for some $x_1\in M^\perp$ and $x_2\in \mathcal{N}(WP_M)$. As $x\in M$, this gives $x=P_Mx=P_Mx_2$, therefore
\begin{align*}
  Wx=WP_Mx_2=0,
\end{align*}
which is combined with \eqref{useful-01} to obtain
$T_{21}x-T_{22}Cx=0$.
The arbitrariness of $x$ in $M$ ensures that $C$ is a solution of \eqref{defn:(S,T)complementable-01}. Similar reasoning shows that $D$ is a solution of \eqref{defn:(S,T)complementable-02}.

Let $M_r\in \mathcal{L}(H)$ and $M_l\in \mathcal{L}(K)$ be defined by
$$M_r=U_{P_{M}}^*\left(
                            \begin{array}{cc}
                              0 & 0 \\
                            C & I_{M^\perp} \\
                            \end{array}
                          \right)U_{P_{M}},
\quad
M_l=U_{P_{N}}^*\left(
                            \begin{array}{cc}
                              0 & D^* \\
                              0 & I_{N^\perp} \\
                            \end{array}
                          \right)U_{P_{N}}.
$$
A simple computation yields the validity of \eqref{equ:M_s and M_T-01} and \eqref{equ:M_s and M_T-02}.
\end{proof}

\begin{remark} Suppose that the operator $T_{22}$ defined by \eqref{equ:A matrix-02} has the polar decomposition. Then from Lemmas~\ref{lem:dougla's theorem} and \ref{lem:polar decomposition of T star}, it is easily seen that conditions stated in item (ii) of Proposition~\ref{prop of S T complementable}
can be simplified as
\begin{equation*}\mathcal{R}(T_{21})\subseteq \mathcal{R}(T_{22})\quad\mbox{and}\quad \mathcal{R}(T_{12}^*)\subseteq \mathcal{R}(T_{22}^*),
\end{equation*}
which is originally characterized in \cite[Proposition~3.2]{ACS} for Hilbert space operators.
\end{remark}

\section{$(M,N)$-weakly complementable operators}\label{sec:weak complementarity}

In this section, we always assume that $M$ and $N$ are orthogonally complemented closed submodules of $H $ and $K$, respectively.

\begin{definition}\label{defn:(S,T)-weakly}
 An operator $T\in \mathcal{L}(H,K)$ is said to be  $(M,N)$-weakly complementable if all the systems
\begin{align}\label{equ:A operator equation-01}
   &V_{T_{22}}Y_{1}=T_{21}, \quad Y_{1}\in\mathcal{L}(M,M^\perp  ),\\
   \label{equ:A operator equation-02}&|T_{22}|^{\frac{1}{2}}Y_{2}=T_{12}^*,\quad Y_{2}\in\mathcal{L}(N,M^\perp  ), \\
   \label{equ:A operator equation-03} &|T_{22}^*|^{\frac{1}{2}}Z_1=T_{21},\quad Z_{1}\in\mathcal{L}(M,N^\perp),\\
    \label{equ:A operator equation-04}&V_{T_{22}}^*Z_{2}=T_{12}^*,\quad Z_{2}\in\mathcal{L}(N,N^\perp  )
\end{align}
have the reduced solutions, denoted by $E,F,\widetilde{E}$ and $\tilde{F}$ respectively,
where $T_{12}, T_{21}$ and $T_{22}$ are defined by \eqref{equ:A matrix-01} and \eqref{equ:A matrix-02}, while $V_{T_{22}}$ is the reduced solution of the system
$$|T_{22}^*|^\frac12 X=T_{22},\quad X\in\mathcal{L}(M^\perp,N^\perp).$$
\end{definition}

\begin{remark}Suppose that $T\in \mathcal{L}(H,K)$ is $(M,N)$-weakly complementable. For simplicity, in the remaining of this section the notations $T_{ij} (1\le i,j\le 2)$, $E,F,\widetilde{E}$ and $\tilde{F}$ will be used without being specified. From Definition~\ref{defn:(S,T)-weakly},
  \eqref{range closure of V T}, \eqref{red sys T++} together with   Lemmas~\ref{lem:Range closure of TT and T and closures are the same-01} and \ref{lem:Range closure of TT and T and closures are the same-02}, we have
\begin{align}\label{red-for citation-01}   &V_{T_{22}}E=T_{21}, \quad \mathcal{R}(E)\subseteq \overline{\mathcal{R}(V_{T_{22}}^*)}=\overline{\mathcal{R}(T_{22}^*)},\\
   \label{red-for citation-02} &|T_{22}|^{\frac{1}{2}}F=T_{12}^*,\quad \mathcal{R}(F)\subseteq \overline{\mathcal{R}(|T_{22}|^{\frac{1}{2}})}=\overline{\mathcal{R}(T_{22}^*)}, \\
   \label{red-for citation-03} &|T_{22}^*|^{\frac{1}{2}}\widetilde{E}=T_{21},\quad \mathcal{R}(\widetilde{E})\subseteq \overline{\mathcal{R}(|T_{22}^*|^{\frac{1}{2}})}=\overline{\mathcal{R}(T_{22})},\\
    \label{red-for citation-04} &V_{T_{22}}^*\widetilde{F}=T_{12}^*,\quad \mathcal{R}(\widetilde{F})\subseteq \overline{\mathcal{R}(V_{T_{22}})}=\overline{\mathcal{R}(T_{22})},\\
     \label{red-for citation-05}&|T^*_{22}|^\frac12 V_{T_{22}}=T_{22}, \quad |T_{22}|^\frac12 V_{T_{22}}^*=T^*_{22}.
\end{align}

\end{remark}

\begin{definition}\label{defn:A/(S,T)}
Suppose that $T\in \mathcal{L}(H,K)$ is $(M,N)$-weakly complementable.  The bilateral shorted operator of $T$ with respect to $(M,N)$, written $T_{/(M,N)}$, is the operator in $\mathcal{L}(H,K)$ given by
\begin{equation}\label{formula with half} T_{/(M,N)}=U^*_{P_N}\left(
               \begin{array}{cc}
                 T_{11}-\frac12\big[F^*E+\widetilde{F}^*\widetilde{E}\big] &\quad 0 \\
                 0 &\quad 0 \\
               \end{array}
             \right)U_{P_M}.
\end{equation}
\end{definition}

A simple combination of  Definitions \ref{defn:(S,T)-weakly} and \ref{defn:A/(S,T)}   with Corollary~\ref{cor:T and Vt} gives  the following proposition.
\begin{proposition}
Suppose that $T\in \mathcal{L}(H,K)$ is $(M,N)$-weakly complementable. Then $T^*$ is $(N,M)$-weakly complementable such that  $\big(T_{/(M,N)}\big)^*=T^*_{/(N,M)}$.
\end{proposition}

\begin{proposition}{\rm \cite[Definitions~3.5 and 4.1]{ACS}}\label{prop:polar decomposition-weak reduced}
Let $T\in \mathcal{L}(H,K)$ be given such that $T_{22}$  has the polar decomposition. Then $T$ is $(M,N)$-weakly complementable if and only if
\begin{equation}\label{range of R(A21) and R(A12*)}
  \mathcal{R}(T_{12}^*)\subseteq \mathcal{R}(|T_{22}|^\frac12), \quad \mathcal{R}(T_{21})\subseteq \mathcal{R}(|T_{22}^*|^\frac12).
\end{equation}
In such case,
\begin{equation}\label{simplified formula without half}T_{/(M,N)}=U^*_{P_N} \left(
               \begin{array}{cc}
                 T_{11}-F^*E &\quad 0 \\
                 0 &\quad 0 \\
               \end{array}
             \right)U_{P_M}.
\end{equation}
\end{proposition}
\begin{proof}Suppose that $T$ is $(M,N)$-weakly complementable. From \eqref{red-for citation-02} and \eqref{red-for citation-03}, it is clear that \eqref{range of R(A21) and R(A12*)} is satisfied.

Conversely, suppose that \eqref{range of R(A21) and R(A12*)} is satisfied. Since $T_{22}$ has the polar decomposition,
by Lemma~\ref{lem:polar decomposition of T star} $\overline{\mathcal{R}(T_{22})}$ and $\overline{\mathcal{R}(T_{22}^*)}$ are
orthogonally complemented in $N^\perp$ and $M^\perp$, respectively. Hence, it can be deduced from \eqref{range of R(A21) and R(A12*)} and Lemma~\ref{lem:dougla's theorem} that
systems \eqref{equ:A operator equation-02} and \eqref{equ:A operator equation-03} have the reduced solutions $F$ and $\widetilde{E}$, respectively.
Let the polar decomposition of $T_{22}$ be represented by
$$T_{22}=U_{22}|T_{22}|,\quad U_{22}^*U_{22}=P_{\overline{\mathcal{R}(T_{22}^*)}},$$
where $U_{22}\in \mathcal{L}(M^\perp,N^\perp)$ is a partial isometry. A direct use of \eqref{V and T-02} yields
\begin{equation*}V_{T_{22}}U_{22}^*=|T_{22}^*|^\frac12,\end{equation*}
which is combined with  \eqref{red-for citation-03} to get
$$V_{T_{22}}(U_{22}^*\widetilde{E})=|T_{22}^*|^\frac12\widetilde{E}=T_{21}.$$
Since $\mathcal{R}(U_{22}^*\widetilde{E})\subseteq \mathcal{R}(U_{22}^*)=\overline{\mathcal{R}(T_{22}^*)}$ and $\overline{\mathcal{R}(T_{22}^*)}=\overline{\mathcal{R}(V^*_{T_{22}})}$ (see \eqref{range closure of V T}), we see that
$U_{22}^*\widetilde{E}$ is the reduced solution of \eqref{equ:A operator equation-01}. Similarly, it can be deduced from the second equation in \eqref{V and T-03} together with
\eqref{red-for citation-02} that $U_{22}F$ is the reduced solution of \eqref{equ:A operator equation-04}. So, $T$ is $(M,N)$-weakly complementable such that \begin{equation}\label{relations E F etc}\widetilde{F}=U_{22}F, \quad E=U_{22}^*\widetilde{E}.\end{equation}
The second equation above together with \eqref{red-for citation-03} yields
$$U_{22}E=U_{22}U_{22}^*\widetilde{E}=\widetilde{E}\quad\mbox{and}\quad U_{22}^*U_{22}E=E,$$
which are combined with the first equation in \eqref{relations E F etc} to get
$$\widetilde{F}^*\widetilde{E}=F^*U_{22}^*U_{22}E=F^*E.$$ So,  formula \eqref{formula with half} can be simplified as \eqref{simplified formula without half}.
\end{proof}

\begin{proposition}\label{prop:s to w-01}
Every $(M,N)$-complementable operator is $(M,N)$-weakly complementable, and its bilateral shorted operator can be expressed as \eqref{simplified formula without half}.
  \end{proposition}
 \begin{proof}Let $T\in\mathcal{L}(H,K)$ be $(M,N)$-complementable such that $U_{P_{N}}TU_{P_{M}}^*$ and
 $ T_{/(M,N)} $
 are given by \eqref{2 form A} and \eqref{formula with half}, respectively. By Proposition~\ref{prop of S T complementable},
 there exist $C\in\mathcal{L}(M,M^\perp)$ and $D\in\mathcal{L}(N,N^\perp)$ such that
 \begin{equation}\label{equs for C and D}T_{22}C=T_{21}\quad\mbox{and}\quad T_{22}^*D=T_{12}^*.\end{equation}
The equations above together with \eqref{red-for citation-05} yield
 \begin{align*}
   &T_{21}=\big(|T_{22}|^\frac{1}{2}V_{T_{22}}^*\big)^* C=V_{T_{22}}\big(|T_{22}|^\frac{1}{2}C\big),\quad T_{12}^*=|T_{22}|^\frac{1}{2}\big(V_{T_{22}}^*D\big),\\
  &T_{21} =|T_{22}^*|^\frac{1}{2}(V_{T_{22}}C),\quad  T_{12}^*=\big(|T_{22}^*|^\frac{1}{2}V_{T_{22}}\big)^*D=V_{T_{22}}^*\big(|T_{22}^*|^\frac{1}{2}D\big).
\end{align*}
Hence,  it can be deduced from \eqref{range closure of V T} that
$$|T_{22}|^\frac{1}{2}C,\quad V_{T_{22}}^*D, \quad V_{T_{22}}C, \quad |T_{22}^*|^\frac{1}{2}D$$ are the reduced solutions of the systems \eqref{equ:A operator equation-01}--\eqref{equ:A operator equation-04}, respectively. Therefore, $T$ is $(M,N)$-weakly complementable.
In virtue of \eqref{V T commutative}, we have
\begin{align}\label{equ:FE}\widetilde{F}^*\widetilde{E}=D^*|T_{22}^*|^\frac{1}{2}V_{T_{22}}C=D^*V_{T_{22}}|T_{22}|^\frac{1}{2}C=F^*E.
\end{align}
So, the desired conclusion follows.
\end{proof}

\begin{proposition}{\rm (cf.\,\cite[Corollary~4.6]{ACS})}\label{cor:range A(M,N)}
Suppose that $T\in\mathcal{L}(H,K)$ is $(M,N)$-weakly complementable. Then
\begin{equation*}
  \mathcal{R}(T)\cap N\subseteq \mathcal{R}\big(T_{/(M,N)}\big)\subseteq\overline{\mathcal{R}(T)}\cap N.
  \end{equation*}
  \end{proposition}
\begin{proof}Let  $U_{P_{N}}TU_{P_{M}}^*$ and
 $T_{/(M,N)}$
 be given by \eqref{2 form A} and \eqref{formula with half}, respectively.

First, we prove that $\mathcal{R}(T)\cap N\subseteq \mathcal{R}\big(T_{/(M,N)}\big)$.
Given an arbitrary element $u$ in $\mathcal{R}(T)\cap N$, we have $$U_{P_{N}}u=(u,0)^T\quad \mbox{and}\quad u=Tz$$ for some $z\in H$, hence
$U_{P_{N}}TU_{P_{M}}^*\cdot U_{P_{M}}z=U_{P_{N}}u$, i.e.,
\begin{equation*}
  \left(
     \begin{array}{cc}
       T_{11} & T_{12} \\
       T_{21} & T_{22} \\
     \end{array}
   \right)
   \left(
     \begin{array}{c}
       x \\
       y \\
     \end{array}
   \right)
   =\left(
      \begin{array}{c}
        u \\
        0 \\
      \end{array}
    \right),
\end{equation*}
where $x=P_{M}(z)$ and $y=(I_H-P_{M})(z)$.
It follows that
\begin{equation}\label{2 equs A 12 21}T_{11}x+T_{12}y=u, \quad T_{21}x+T_{22}y=0.\end{equation}
The second equation above, together with \eqref{red-for citation-01} and the second equation in \eqref{red-for citation-05}, yields
$$0=V_{T_{22}}Ex+V_{T_{22}}|T_{22}|^\frac12 y=V_{T_{22}}w,$$
where
$$w=Ex+|T_{22}|^{\frac{1}{2}}y\in \mathcal{R}(E)+\mathcal{R}\big(|T_{22}|^{\frac{1}{2}}\big)\subseteq\overline{\mathcal{R}(V_{T_{22}}^*)}.$$
As a result, $$w\in\mathcal{N}(V_{T_{22}})\cap \overline{\mathcal{R}(V_{T_{22}}^*)}=\{0\}.$$
Consequently, $|T_{22}|^{\frac{1}{2}}y=-Ex$ and thus we may use \eqref{red-for citation-02} to obtain
\begin{align*}T_{12}y=F^*|T_{22}|^{\frac{1}{2}}y=-F^*Ex.\end{align*}
So the first equation in \eqref{2 equs A 12 21} turns out to be
\begin{equation}\label{1st form of u}
   u=T_{11}x-F^*Ex.
\end{equation}
Here the point is, $u$ can be expressed alternatively as
\begin{equation}\label{2nd form of u}u=T_{11}x-\widetilde{F}^*\widetilde{E}x.
\end{equation}
Indeed, by \eqref{red-for citation-03} and the first equation in \eqref{red-for citation-05}, it is easily seen that the second equation in \eqref{2 equs A 12 21}
can be rephrased as
$$|T_{22}^*|^{\frac{1}{2}}\big(\widetilde{E}x+V_{T_{22}}y\big)=0,$$
which, as shown in the derivation of $w=0$, implies that $\widetilde{E}x+V_{T_{22}}y=0$. This, together with the first equation in \eqref{2 equs A 12 21} and \eqref{red-for citation-04}, yields
\begin{equation*}
   u=T_{11}x+\widetilde{F}^*V_{T_{22}}y=T_{11}x-\widetilde{F}^*\widetilde{E}x.
\end{equation*}
Combining \eqref{1st form of u} and \eqref{2nd form of u}  gives
$$u=\Big[T_{11}-\frac12\big(F^*E+\widetilde{F}^*\widetilde{E}\big)\Big]x.$$
It follows that
\begin{align*}u=&U_{P_N}^*\left(
               \begin{array}{c}
                 u \\
                 0 \\
               \end{array}
             \right)=U_{P_N}^*\left(
               \begin{array}{cc}
                 T_{11}-\frac12\big[F^*E+\widetilde{F}^*\widetilde{E}\big] &\quad 0 \\
                 0 &\quad 0 \\
               \end{array}
             \right)\left(
                      \begin{array}{c}
                        x \\
                        0 \\
                      \end{array}
                    \right)\\
=&T_{/(M,N)}\cdot U_{P_M}^*\left(
                                                 \begin{array}{c}
                                                   x \\
                                                   0 \\
                                                 \end{array}
                                               \right)\in \mathcal{R}\big(T_{/(M,N)}\big).
\end{align*}
This shows the validity of $\mathcal{R}(T)\cap N\subseteq \mathcal{R}\big(T_{/(M,N)}\big)$.

Next, we prove  that $\mathcal{R}(T_{/(M,N)})\subseteq\overline{\mathcal{R}(T)}\cap N$. From \eqref{formula with half}, it is clear that
$$\mathcal{R}(T_{/(M,N)})\subseteq N,\quad T_{/(M,N)}=\frac12 (S+\widetilde{S}),$$ where $S,\widetilde{S}\in \mathcal{L}(H,K)$  are given by
 $$S=U^*_{P_N}\left(
    \begin{array}{cc}
     T_{11}-F^*E &\quad 0 \\
      0 &\quad 0 \\
    \end{array}
  \right)U_{P_M}, \quad
\widetilde{S}=U^*_{P_N}\left(
    \begin{array}{cc}
      T_{11}-\widetilde{F}^*\widetilde{E} &\quad 0 \\
      0 &\quad 0 \\
    \end{array}
  \right)U_{P_M}.$$
So, it is sufficient to verify that both
$\mathcal{R}(S)$ and $\mathcal{R}(\widetilde{S})$ are  contained in $\overline{\mathcal{R}(T)}$.

For each $x\in M$, we have $-Ex\in \mathcal{R}(E)$, so by \eqref{red-for citation-01}
there exists a sequence $\{y_n\}$ in $M^\bot$ such that
\begin{equation*}\label{equ:limit y_n}
  -Ex=\lim_{n\to \infty}|T_{22}|^\frac12y_n,
\end{equation*}
which leads  by \eqref{red-for citation-01} and \eqref{red-for citation-05} to
$$\lim_{n\to \infty}(T_{21}x+T_{22}y_n)=\lim_{n\to \infty}V_{T_{22}}\big(Ex+|T_{22}|^\frac12y_n\big)=0.$$
So for any $x\in M$ and $y\in M^\perp$, we have
\begin{align*}U_{P_N}SU^*_{P_M} \left(
                                  \begin{array}{c}
                                    x \\
                                    y \\
                                  \end{array}
                                \right)
=&\begin{pmatrix}
                                           T_{11}x-F^*Ex \\
                                           0 \\
                                         \end{pmatrix}=\lim_{n\to\infty}\begin{pmatrix}
                                           T_{11}x+F^*|T_{22}|^\frac12 y_n \\
                                           0 \\
                                          \end{pmatrix}\\
=&\lim_{n\to\infty}\begin{pmatrix}
                                           T_{11}x+T_{12}y_n \\
                                          T_{21}x+T_{22}y_n \\
                                         \end{pmatrix}=\lim_{n\to\infty}U_{P_N}TU^*_{P_M}
                                 \left(
                                   \begin{array}{c}
                                     x \\
                                     y_n \\
                                   \end{array}
                                 \right).
\end{align*}
As $U_{P_N}$ is a unitary, this shows that $\mathcal{R}(S)\subseteq \overline{\mathcal{R}(T)}$.

Similarly, for each $x\in M$ we have
$$-\widetilde{E}x\in \mathcal{R}(\widetilde{E})\subseteq \overline{\mathcal{R}(|T_{22}^*|^\frac12)} = \overline{\mathcal{R}(V_{T_{22}})},$$
so there exists  $\{t_n\}$ in $M^\bot$ such that
$-\widetilde{E}x=\lim\limits_{n\to \infty}V_{T_{22}}t_n$.
Following the same line in the derivation of
$\mathcal{R}(S)\subseteq \overline{\mathcal{R}(T)}$, we may use \eqref{red-for citation-04}, \eqref{red-for citation-03} and the first equation in \eqref{red-for citation-05} to conclude that
\begin{align*}U_{P_N}\widetilde{S}U^*_{P_M}\left(
                                             \begin{array}{c}
                                               x \\
                                               y \\
                                             \end{array}
                                           \right)
=\lim_{n\to\infty}U_{P_N}TU^*_{P_M}
                                 \left(
                                   \begin{array}{c}
                                     x \\
                                     t_n \\
                                   \end{array}
                                 \right)
\end{align*}
for any $x\in M$ and $y\in M^\perp$. Hence,  $\mathcal{R}(\widetilde{S})\subseteq \overline{\mathcal{R}(T)}$ as desired.
\end{proof}

\begin{corollary}\label{cor:range characterization00}{\rm (cf.\,\cite[Corollary~4.9]{ACS})}
Suppose that $T\in \mathcal{L}(H,K)$ is $(M, N)$-complementable. Then
\begin{equation}\label{range-kernel-01}\mathcal{R}\big(T_{/(M,N)}\big)=\mathcal{R}(T)\cap N, \quad \mathcal{N}\big(T_{/(M,N)}\big)=M^\perp+ \mathcal{N}(T).\end{equation}
\end{corollary}
\begin{proof} By Proposition~\ref{prop of S T complementable}(ii),
 there exist $C\in\mathcal{L}(M,M^\perp)$ and $D\in\mathcal{L}(N,N^\perp)$ such that
 \eqref{equs for C and D} is satisfied.

First, we prove the first equation in \eqref{range-kernel-01}. By Proposition~\ref{cor:range A(M,N)} and the observation $\mathcal{R}\big(T_{/(M,N)}\big)\subseteq N$,
it needs only to verify that $\mathcal{R}\big(T_{/(M,N)}\big)\subseteq \mathcal{R}(T)$. Given $z \in H$, let
$x=P_{M}(z)$ and $y=(I_H-P_{M})(z)$. That is, $U_{P_M}z=(x,y)^T$. By Proposition~\ref{prop:s to w-01},
\begin{equation}\label{a.e.0}
  U_{P_{N}}T_{/(M,N)}U_{P_{M}}^*\cdot U_{P_{M}}z                                               =\left(
                                                  \begin{array}{c}
                                                    T_{11}x-F^*Ex \\
                                                    0 \\
                                                  \end{array}
                                                \right).
\end{equation}
Let
 $u\in M$, $v\in M^\perp$ and $w\in H$ be defined by
 $$u=x, \quad v=-Cx,\quad w=u+v.$$ It can be deduced from \eqref{equs for C and D},  \eqref{red-for citation-05} and  \eqref{equ:FE} that
\begin{align*}
  &T_{11}u+T_{12}v=T_{11}x-T_{12}Cx=T_{11}x-D^*V_{T_{22}}|T_{22}|^\frac{1}{2}Cx=T_{11}x-F^*Ex,\\
  &T_{21}u+T_{22}v=T_{22}Cu+T_{22}v=T_{22}Cx-T_{22}Cx=0.
\end{align*}
This shows that for any $z\in H$, there exists $w\in H$ such that
\begin{equation*}
U_{P_{N}}TU_{P_{M}}^*\cdot U_{P_{M}}w=U_{P_{N}}T_{/(M,N)}U_{P_{M}}^*\cdot U_{P_{M}}z.
\end{equation*}
 By the arbitrariness of $z$, we obtain $\mathcal{R}(T_{/(M,N)})\subseteq \mathcal{R}(T)$.

Next, we prove the second equation in \eqref{range-kernel-01}. By \eqref{simplified formula without half}, it is clear that $M^\perp\subseteq\mathcal{N}\big(T_{/(M,N)}\big).$ For each $z\in \mathcal{N}(T)$, let
$x=P_{M}(z)$ and $y=(I_H-P_{M})(z)$ as before. Exploring $U_{P_N}TU_{P_M}^*\cdot U_{P_M}z=0$ yields
\begin{align*}
  &0=T_{11}x+T_{12}y=T_{11}x+D^*T_{22}y=T_{11}x+D^*|T_{22}^*|^\frac{1}{2}V_{T_{22}}y,\\
  &0=T_{21}x+T_{22}y=T_{22}Cx+T_{22}y=|T_{22}^*|^\frac12V_{T_{22}}(Cx+y).
\end{align*}
Hence
$$V_{T_{22}}(Cx+y)\in \overline{\mathcal{R}(V_{T_{22}})}\cap  \mathcal{N}(|T_{22}^*|^\frac12)=\overline{\mathcal{R}(T_{22})}\cap \mathcal{N}(T_{22}^*)=\{0\}.$$
Therefore, $V_{T_{22}}y=-V_{T_{22}}Cx$, and thus by \eqref{equ:FE} we have
$$0=T_{11}x+T_{12}y=T_{11}x-D^*|T_{22}^*|^\frac{1}{2}V_{T_{22}}Cx=T_{11}x-F^*Ex,$$
which leads by \eqref{a.e.0} to
$$U_{P_{N}}T_{/(M,N)}z=U_{P_{N}}T_{/(M,N)}U_{P_{M}}^*\cdot U_{P_{M}}z=(0,0)^T.$$
So, $\mathcal{N}(T)\subseteq\mathcal{N}\big(T_{/(M,N)}\big)$. This shows that $M^\perp+\mathcal{N}(T)\subseteq\mathcal{N}\big(T_{/(M,N)}\big).$

Suppose now that $w\in \mathcal{N}\big(T_{/(M,N)}\big)$. Let  $$u=P_{M}(w),\quad v=(I_H-P_M)(w).$$ Then from the above arguments, we have
\begin{align*}0=&T_{11}u-F^*Eu=T_{11}u-D^*|T_{22}^*|^\frac{1}{2}V_{T_{22}}Cu\\
=&T_{11}u-D^*T_{22}Cu=T_{11}u+T_{12}(-Cu),\\
0=&T_{22}Cu+T_{22}(-Cu)=T_{21}u+T_{22}(-Cu).
\end{align*}
Therefore,
\begin{equation*}
  \left(
     \begin{array}{cc}
       T_{11} & T_{12} \\
       T_{21} & T_{22} \\
     \end{array}
   \right)\left(
            \begin{array}{c}
              u \\
              -Cu \\
            \end{array}
          \right)=\left(
                    \begin{array}{c}
                      0 \\
                      0 \\
                    \end{array}
                  \right),
\end{equation*}
which means that $u-Cu\in\mathcal{N}(T)$. Hence,
$$w=u+v=(u-Cu)+(Cu+v)\in \mathcal{N}(T)+M^\perp.$$
This shows that $\mathcal{N}\big(T_{/(M,N)}\big)\subseteq\mathcal{N}(T)+M^\perp$. So, the desired equation follows.
\end{proof}

\section{Some investigations around two special positive operators}

\subsection{Solvability of two special operator equations}\label{subsec:two positive operators}
In this subsection, we will construct two special positive operators $A_0$ and $B_0$, and study the  solvability of equations \eqref{equ:operator equation of(A+B)X=B} and
\eqref{equ:half operator equation}.
  Let $\ell^2(\mathbb{N})$ be the Hilbert space consisting of  all complex sequence $x=(x_1,x_2,\cdots,x_n,\cdots)$ such that $\|x\|^2=\sum\limits_{i=1}^{\infty}|x_i|^2<\infty$,  and let $S\in \mathbb{B}\big(\ell^2(\mathbb{N})\big)$ be a diagonal operator defined by
\begin{equation}\label{operator s for A B}S(x)=\left(x_1,\frac{1}{2}x_2,\cdots,\frac{1}{n}x_n,\cdots\right)\end{equation} for each $x=(x_1,x_2,\cdots,x_n,\cdots)\in \ell^2(\mathbb{N})$. It is clear that
$$\|S\|=1,\quad S\geq 0, \quad \mathcal{N}(S)=0.$$ Put $H=\ell^2(\mathbb{N})\oplus \ell^2(\mathbb{N})$ and define $A_0,B_0\in\mathbb{B}(H)$ by setting
\begin{equation}\label{counterexample A B}A_0=\left(
                       \begin{array}{cc}
                         I & S \\
                         S & S^2 \\
                       \end{array}
                     \right),\quad B_0=\left(
      \begin{array}{cc}
        I & 0 \\
        0 & 0 \\
      \end{array}
    \right),\end{equation}
where $I$ denotes the identity operator on $\ell^2(\mathbb{N})$.

\begin{proposition}\label{prop:R B not contained} The operators $A_0$ and $B_0$ defined by \eqref{counterexample A B}
  are both positive, and
the  equation
\begin{equation}\label{equ:operator equation of(A+B)X=B}
  (A_0+B_0)X=B_0,\quad X\in\mathbb{B}(H)
\end{equation}
has no solution.
\end{proposition}
\begin{proof} The positivity of $B_0$ is obvious. Rewrite $A_0$ as
$$A_0=\left(
        \begin{array}{cc}
          I & S \\
        \end{array}
      \right)^*
\left(
        \begin{array}{cc}
          I & S \\
        \end{array}
      \right),$$
the positivity of $A_0$ follows. To show the insolvability of \eqref{equ:operator equation of(A+B)X=B}, it needs only to prove that $\mathcal{R}(B_0)\nsubseteq \mathcal{R}(A_0+B_0)$.
Evidently, $\mathcal{R}(B_0)=\ell^2(\mathbb{N})\oplus \{0\}$. Specifically,  $(x_0,0)^T\in\mathcal{R}(B_0)$, where
$x_0\in\ell^2(\mathbb{N})$ is given by
 $x_0=(1,\frac{1}{2},\cdots,\frac{1}{n},\cdots)$.

 We claim that $(x_0,0)^T\notin\mathcal{R}(A_0+B_0)$. Suppose on the contrary that there exist
$\xi,\eta\in \ell^2(\mathbb{N})$ such that $(A_0+B_0)(\xi,\eta)^T=(x_0,0)$.
Then $$2\xi+S\eta=x_0,\quad S(\xi+S\eta)=0,$$
since by \eqref{counterexample A B} we have
\begin{equation}\label{patitioned A plus B}A_0+B_0=\left(
                       \begin{array}{cc}
                         2I & S \\
                         S & S^2 \\
                       \end{array}
                     \right).\end{equation}
As $S$ is injective, it gives $\xi+S\eta=0$ and thus
$$x_0=2\xi+S\eta=-S\eta.$$ It follows from \eqref{operator s for A B} that $\eta=(-1,-1,\cdots,-1,\cdots) \in \ell^2(\mathbb{N})$, which is a  contradiction. This shows that
$\mathcal{R}(B_0)\nsubseteq \mathcal{R}(A_0+B_0)$.
\end{proof}

\begin{proposition}Let $A_0,B_0\in\mathbb{B}(H)$ be defined  by \eqref{counterexample A B}.    Then
\begin{equation}\label{equ:(A+B)half}(A_0+B_0)^{\frac12}=\left(
                                  \begin{array}{cc}
                                    f(S)^{-1}(S+2I) & f(S)^{-1}S \\
                                    f(S)^{-1}S&f(S)^{-1}S(S+I) \\
                                  \end{array}
                                \right),
\end{equation}
where  $S$ is defined by \eqref{operator s for A B}, and
\begin{equation}\label{equ:defn of f t}f(t)=\sqrt{t^2+2t+2},\quad t\in [0,+\infty).\end{equation}
\end{proposition}
\begin{proof}Denote by $\mathbb{R}^{2\times 2}$ the set of all $2\times 2$ real  matrices.
For each $t\ge 0$, let $R_t\in \mathbb{R}^{2\times 2}$ be defined by
$$R_t=\left(
        \begin{array}{cc}
          2 &  t\\
         t  & t^2 \\
        \end{array}
      \right).$$
Obviously,  $R_t$ is positive semi-definite.  It is routine to verify that
$$R_t^\frac12=\frac{1}{f(t)}\left(
        \begin{array}{cc}
          t+2 &  t\\
         t  & t(t+1)\\
        \end{array}
      \right),\quad  t\in [0,+\infty).$$
Since $S$ is positive and $A_0+B_0$ is partitioned as \eqref{patitioned A plus B}, the desired conclusion follows.
\end{proof}

\begin{proposition}\label{prop:solution of (A+B)12X=B}Let $A_0,B_0\in\mathbb{B}(H)$ be defined by \eqref{counterexample A B}.
Then the equation
\begin{equation}\label{equ:half operator equation}
  (A_0+B_0)^\frac12 X=B_0,\quad X\in\mathbb{B}(H)
\end{equation}
has a unique solution formulated by
\begin{equation}\label{unique solu X}X=\left(
                    \begin{array}{cc}
                      f(S)^{-1}(S+I) &\quad 0 \\
                      -f(S)^{-1} &\quad 0 \\
                    \end{array}
                  \right),
\end{equation}
where  $S$ and $f$ are defined by \eqref{operator s for A B} and \eqref{equ:defn of f t}, respectively.
\end{proposition}
\begin{proof}Suppose that $X=(X_{ij})\in \mathbb{B}(H)$ is a solution of \eqref{equ:half operator equation}, in which $X_{ij}\in \mathbb{B}\big(\ell^2(\mathbb{N})\big)$. Since $f(S)$ is invertible, by \eqref{equ:(A+B)half}, \eqref{equ:half operator equation} and \eqref{counterexample A B} we have
\begin{align*}&(S+2I)X_{11}+SX_{21}=f(S),\\
&(S+2I)X_{12}+SX_{22}=0,\\
&S\left[X_{11}+(S+I)X_{21}\right]=0,\\
&S\left[X_{12}+(S+I)X_{22}\right]=0.
\end{align*}
Utilizing the injectivity of $S$ and the invertibility of $f(S)$, it is easy to obtain
\begin{align*}&X_{22}=X_{12}=0,\quad X_{21}=-f(S)^{-1},\quad  X_{11}=f(S)^{-1}(S+I).
\end{align*}
Hence, the conclusion follows.
\end{proof}

\subsection{A remark on the  weak complementarity and  complementarity}\label{subsec:complementarity}

The purpose of this subsection is to show the inconsistency of the weak complementarity and  complementarity.
Given Hilbert spaces $H_1$ and $H_2$, let $H=H_1\oplus H_2$, $M=H_1\oplus \{0\}$
and $P_M: H\to M$ be the projection from $H$ onto $M$. It is obvious that $M^\perp={0}\oplus H_2$. Define the unitaries $\pi_1:H_1\rightarrow M$ and $\pi_2:H_2\rightarrow M^\perp $ by setting
\begin{equation*}
  \pi_1(x)=\left(
             \begin{array}{c}
               x \\
               0 \\
             \end{array}
           \right),\quad
\pi_2(y)=\left(
             \begin{array}{c}
               0 \\
               y \\
             \end{array}
           \right)
\quad (x\in H_1,y\in H_2).
\end{equation*}
Let $T\in \mathbb{B}(H)$ be an arbitrary operator partitioned as
\begin{equation*}\label{W 2times2}
  T=\left(
                          \begin{array}{cc}
                            T_{11} & T_{12} \\
                            T_{21} & T_{22} \\
                          \end{array}
                        \right),
\end{equation*}
in which $T_{ij}\in \mathbb{B}(H_j,H_i)$ ($1\le i,j\le 2$).
It follows  easily from \eqref{2 form A}--\eqref{equ:A matrix-02} that
\begin{equation*}
  U_{P_{M}}TU_{P_{M}}^*=\left(
                          \begin{array}{cc}
                            \widetilde{T_{11}} & \widetilde{T_{12}} \\
                            \widetilde{T_{21}} & \widetilde{T_{22}} \\
                          \end{array}
                        \right),
\end{equation*}where \begin{align}
  \label{W11 and W12}&\widetilde{T_{11}}=\pi_1T_{11}\pi_1^*\in \mathbb{B}(M), \quad \widetilde{T_{12}}=\pi_1T_{12}\pi_2^*\in \mathbb{B}(M^\perp,M),  \\
  \label{W21 and W22}&\widetilde{T_{21}}=\pi_2T_{21}\pi_1^*\in \mathbb{B}(M,M^\perp),\quad \widetilde{T_{22}}=\pi_2T_{22}\pi_2^*\in \mathbb{B}(M^\perp).
\end{align}
By Proposition~\ref{prop of S T complementable}, $T$ is $(M,M)$-complementable if and only if each of the systems
\begin{align*}
  & \widetilde{T_{22}}X=\widetilde{T_{21}}, \quad X\in \mathbb{B}(M,M^\perp),\\
  & \widetilde{T_{22}}^*Y=\widetilde{T_{12}}^*, \quad Y\in \mathbb{B}(M,M^\perp)
\end{align*}
has a solution. In virtue of \eqref{W11 and W12} and \eqref{W21 and W22}, the latter condition can be rephrased as
the solvability of the systems
\begin{align*}
  & T_{22}X^\prime=T_{21}, \quad X^\prime\in \mathbb{B}(H_1,H_2),\\
  & T_{22}^*Y^\prime=T_{12}^*, \quad Y^\prime\in \mathbb{B}(H_1,H_2).
\end{align*}
Similarly, when dealing with the $(M,M)$-weakly complementable operators, we need only consider operators $T_{ij}$ rather than operators $\widetilde{T_{ij}}$. Such an observation will be used in the proof of our next theorem without being specified.

\begin{theorem}There exist a closed subspace  $M$ of certain Hilbert space $K$ and  an operator $T\in \mathbb{B}(K)$ such that $T$ is $(M,M)$-weakly complementable, whereas $T$ fails to be $(M,M)$-complementable.
\end{theorem}
\begin{proof} Let $A_0,B_0\in\mathbb{B}(H)$ be defined by \eqref{counterexample A B}. Put $K=H\oplus H$, $M=H\oplus \{0\}$ and define $T=(T_{ij})\in \mathbb{B}(K)$ by setting
$$T_{11}=T_{12}=T_{21}=B_0,\quad T_{22}=A_0+B_0.$$
By Proposition~\ref{prop:solution of (A+B)12X=B}, we know that \eqref{range of R(A21) and R(A12*)} is satisfied. Hence, it can be concluded by
Proposition~\ref{prop:polar decomposition-weak reduced} that $T$ is $(M,M)$-weakly complementable, since every element in $\mathbb{B}(H)$ has the polar decomposition.

 Observe that the conclusion of  Proposition~\ref{prop:R B not contained} can be rephrased as the insolvability of the equation $T_{22}X=T_{21}$ in $\mathbb{B}(H)$.
 Hence by Proposition~\ref{prop of S T complementable}, $T$ is not $(M,M)$-complementable.
 \end{proof}

\subsection{A remark on the existence of the reduced solution}\label{subsec:reduced solution}
When one deals with Hilbert space operators, it can be deduced directly from the Douglas' theorem that an operator equation
has a solution if and only if it has the reduced solution. Our next theorem shows that the same is not true for Hilbert $C^*$-module operators.

\begin{theorem}\label{thm:solvability does not imp reduced solu}
There exists an operator equation in the setting of adjointable operators on Hilbert $C^*$-modules, which has a unique solution,  whereas this unique solution fails to be the reduced solution.
\end{theorem}
\begin{proof}Given an arbitrary Hilbert space $H$, we may consider $\mathbb{B}(H)$  as a Hilbert $C^*$-module over itself via the inner product
$\langle X,Y\rangle=X^*Y$ for $X,Y\in \mathbb{B}(H)$. Since $\mathbb{B}(H)$ has a unit, we can  identity $\mathcal{L}\big(\mathbb{B}(H)\big)$ with $\mathbb{B}(H)$ in the natural way. To avoid ambiguity, we use the notation $\mathcal{R}_{\mathbb{B}(H)}(T)$
for each $T\in \mathbb{B}(H)$ when $T$ is regarded as an element of $\mathcal{L}\big(\mathbb{B}(H)\big)$.  That is,
$$\mathcal{R}_{\mathbb{B}(H)}(T)=\big\{TY:Y\in \mathbb{B}(H)\big\}.$$

Now, let $H$ be specified as $H=\ell^2(\mathbb{N})\oplus \ell^2(\mathbb{N})$, and consider the operator equation \eqref{equ:half operator equation}, where
$A_0,B_0\in\mathbb{B}(H)$ are given by \eqref{counterexample A B}, in which $S$ is defined by \eqref{operator s for A B}.
By Proposition~\ref{prop:solution of (A+B)12X=B}, this equation has a unique solution $X$ formulated by \eqref{unique solu X}.
To show that $X$ is not the reduced solution, it needs to verify that $\mathcal{R}_{\mathbb{B}(H)}(X)\nsubseteq \overline{\mathcal{R}_{\mathbb{B}(H)}\left[(A_0+B_0)^\frac12\right]}$.

Suppose on the contrary that $\mathcal{R}_{\mathbb{B}(H)}(X)\subseteq \overline{\mathcal{R}_{\mathbb{B}(H)}\left[(A_0+B_0)^\frac12\right]}$. Then
we have $X\in \overline{\mathcal{R}_{\mathbb{B}(H)}\left[(A_0+B_0)^\frac12\right]}$, since $X=X\cdot I_H\in \mathcal{R}_{\mathbb{B}(H)}(X)$. So, there exists a sequence $\{W_n\}_{n=1}^\infty$ in $\mathbb{B}(H)$ with
\begin{equation*}W_n=\left(
           \begin{array}{cc}
             X_n & U_n \\
             Y_n & V_n \\
           \end{array}
         \right)\end{equation*}
such that \begin{equation*}\label{lim W=(A+B)12Wn}
  \lim_{n\to \infty}(A_0+B_0)^\frac12W_n=X
\end{equation*}
in the norm topology. Substituting \eqref{equ:(A+B)half} and \eqref{unique solu X} into the above limitation yields
\begin{align*}
   & (S+2I)X_n+SY_n\to S+I, \nonumber\\
   & (S+2I)U_n+SV_n\to 0,\nonumber\\
   &  SX_n+S(S+I)Y_n\to -I.\\
   & SU_n+S(S+I)V_n\to 0. \nonumber
\end{align*}
By \eqref{operator s for A B}, it is clear that $S$ is a compact operator on $\ell^2(\mathbb{N})$. So from the third limitation as above, it can be deduced that
$I$, the identity operator on $\ell^2(\mathbb{N})$, is a compact operator, which apparently is a contradiction.
\end{proof}

\subsection{Solvability  of equation \eqref{equality is obtained}}\label{subsec:equation wrt parallel sum}
Throughout this subsection, $H$ is a Hilbert space. We use the abbreviations ``sot" and ``wot" to represent ``strong operator topology" and
 ``weak operator topology", respectively. For each pair of  $A$ and $B$ in $\mathbb{B}(H)_+$, it is known (see e.g.\,\cite{ACS}) that their parallel sum  $A:B$ can be clarified  as
$$\left(
    \begin{array}{cc}
      A:B &\quad 0 \\
      0 &\quad 0 \\
    \end{array}
  \right)=\left(
            \begin{array}{cc}
              A &\quad A \\
              A &\quad A+B \\
            \end{array}
          \right)_{/(H\oplus\{0\},H\oplus\{0\})}.
$$
In the special case that $A$ and $B$ are both positive definite, it is well-known that $A:B=A(A+B)^{-1}B=B(A+B)^{-1}A$.
In what follows, we will focus on the solvability of equation \eqref{equality is obtained}.

We begin with an elementary result as follows.
\begin{lemma}\label{lem:some upper bound of (i+x)-1}
Let $X\in \mathbb{B}(H)$ be a positive definite operator. Then
\begin{equation}\label{equ:some upper bound of (i+x)-1}
(I+X)^{-1}\le Y^*Y+(I-Y)^*X^{-1}(I-Y), \quad \forall\, Y \in \mathbb{B}(H).
\end{equation}
Moreover, the equality above occurs if and only if $Y=(I+X)^{-1}$.
\end{lemma}
\begin{proof}By assumption $X$ is positive definite, so $X$ and $I+X$ are both invertible  in  $\mathbb{B}(H)$.
 Let
$$T_{X,Y}=(I+X^{-1})^{\frac{1}{2}}Y-X^{-1}(I+X^{-1})^{-\frac{1}{2}}$$
for each $Y\in \mathbb{B}(H)$.
Then
\begin{align*}T_{X,Y}^*T_{X,Y}=&Y^*(I+X^{-1})Y-Y^*X^{-1}-X^{-1}Y+X^{-2}(I+X^{-1})^{-1}\\
=&Y^*Y+Y^*X^{-1}Y-Y^*X^{-1}-X^{-1}Y+\left[X(X+I)\right]^{-1}\\
=&Y^*Y+(I-Y)^*X^{-1}(I-Y)-X^{-1}+\left[X(X+I)\right]^{-1}\\
=&Y^*Y+(I-Y)^*X^{-1}(I-Y)-(I+X)^{-1}.
\end{align*}
Therefore, \eqref{equ:some upper bound of (i+x)-1} is satisfied and it reduces to be an equality if and only if $T_{X,Y}=0$.
From the definition of $T_{X,Y}$, it is easily seen that
$$T_{X, Y}=0  \Longleftrightarrow Y=X^{-1}(I+X^{-1})^{-1}=(I+X)^{-1}.$$
So, the desired conclusion follows.
\end{proof}

In the special case that $A$ and $B$ are both positive definite, the solvability  of \eqref{equality is obtained} reads as follows.
\begin{lemma}\label{thm:1st on parallel sum}
Suppose that $A,B\in \mathbb{B}(H)$ are both positive definite. Then \eqref{equ: A:B geq C^*AC+(I-C)*B(I-C)}
  is valid and $(A+B)^{-1}B$ is the unique solution of \eqref{equality is obtained}.
\end{lemma}
\begin{proof}Since $A$ and $B$ are both positive definite, we have
\begin{equation}\label{equ:expression of B}A:B=(A^{-1}+B^{-1})^{-1}\quad\mbox{and}\quad B=A^{\frac{1}{2}}X^{-1}A^{\frac{1}{2}},\end{equation} where $X$ is defined by
$$X=A^{\frac{1}{2}}B^{-1}A^{\frac{1}{2}}.$$
Evidently,
\begin{equation*}
  A^{-1}+B^{-1}=A^{-1}+A^{-\frac{1}{2}}XA^{-\frac{1}{2}}=A^{-\frac{1}{2}}(I+X)A^{-\frac{1}{2}},
\end{equation*}
which leads by the first equation in  \eqref{equ:expression of B} to
\begin{equation}\label{equ:A:B=A12(I+X)-1B12}
  A:B=A^{\frac{1}{2}}(I+X)^{-1}A^{\frac{1}{2}}.
\end{equation}
For each $C\in \mathbb{B}(H)$, let
$$Y=A^{\frac{1}{2}}CA^{-\frac{1}{2}}.$$
It is clear that
$$I-Y=A^{\frac{1}{2}}(I-C)A^{-\frac{1}{2}}.$$
The above expressions of $Y$ and $I-Y$ together with the second equation in \eqref{equ:expression of B} yield
\begin{equation*}
  Y^*Y+(I-Y)^*X^{-1}(I-Y)=A^{-\frac{1}{2}}\Big[C^*AC+(I-C)^*B(I-C)\Big]A^{-\frac{1}{2}}.
\end{equation*}
Meanwhile, from \eqref{equ:A:B=A12(I+X)-1B12} we have
$$(I+X)^{-1}= A^{-\frac{1}{2}}(A:B)A^{-\frac{1}{2}}.$$
Hence, a simple use of \eqref{equ:some upper bound of (i+x)-1} gives \eqref{equ: A:B geq C^*AC+(I-C)*B(I-C)}.
Moreover, by Lemma~\ref{lem:some upper bound of (i+x)-1} we know that \eqref{equality is obtained} is true  if and only if
\begin{equation}\label{Y inverse}Y(I+X)=(I+X)Y=I.\end{equation}
From the definitions of $X$ and $Y$,  we see that
\eqref{Y inverse} can be rephrased as
$$C(I+B^{-1}A)=(I+B^{-1}A)C=I.$$
Namely, $$C=(I+B^{-1}A)^{-1}=(A+B)^{-1}B.$$ Thus, the desired conclusion follows.
\end{proof}

\begin{lemma}\label{lem:SOT topology of the projection of the range-1} Suppose that $T\in\mathbb{B}(H)_+$. Let
$T_n=\big(T+\frac{1}{n}I\big)^{-1}T$  for each  $n\in \mathbb{N}$. Then
$\lim\limits_{n\to\infty} T_n=P_{\overline{\mathcal{R}(T)}}$ in sot.
\end{lemma}
\begin{proof}Since every closed subspace of a Hilbert space is always orthogonally complemented, the conclusion is immediate from
\cite[Proposition~2.5]{LLX}.
\end{proof}

Our next theorem gives a generalization of Lemma~\ref{thm:1st on parallel sum}.
\begin{theorem}
Let $A,B\in \mathbb{B}(H)_+$ be such that $\mathcal{R}(B)\subseteq \mathcal{R}(A+B)$. Then
the reduced solution of $(A+B)X=B$ is a solution of \eqref{equality is obtained}.
\end{theorem}
\begin{proof}Let $E\in \mathbb{B}(H)$ be the reduced solution of $(A+B)X=B$. Namely,
$$(A+B)E=B\quad\mbox{and}\quad \mathcal{R}(E)\subseteq \overline{\mathcal{R}(A+B)}.$$
For each $n,m\in\mathbb{N}$, let
\begin{equation*}A_n=A+\frac1n I,\quad B_m=B+\frac1m I.\end{equation*}
By Lemma~\ref{thm:1st on parallel sum},  the following equation
 \begin{equation*}\label{equality is obtained-n m}
  A_n:B_m=X^*A_nX+(I-X)^*B_m(I-X), \quad X\in\mathbb{B}(H)
\end{equation*}
has a unique solution $C_{m,n}$, which can be formulated by
$$C_{m,n}=(A_n+B_m)^{-1}B_m.$$
Taking the limit in sot as $m\to\infty$ gives
\begin{equation}\label{inequality having only n}
  A_n:B=C_n^* A_n C_n+(I-C_n)^*B(I-C_n),
\end{equation}
in which
\begin{align*}C_n=&(A_n+B)^{-1}B=(A_n+B)^{-1}(A+B)E.
\end{align*}
Employing Lemma~\ref{lem:SOT topology of the projection of the range-1}, we arrive at
$$\lim_{n\to\infty} C_n=P_{\overline{\mathcal{R}(A+B)}}E=E\quad\mbox{in sot}.$$
Hence, it follows directly from  \eqref{inequality having only n} that
$$A:B=E^*AE+(I-E)^*B(I-E)$$ by taking the limit in wot. This shows that $E$ is a solution of \eqref{equality is obtained}.
\end{proof}

\begin{remark}Let $A,B\in \mathbb{B}(H)_+$. It is evident that
$$\mathcal{R}(B)\subseteq \mathcal{R}(A+B)\Longleftrightarrow
\mathcal{R}(A)\subseteq \mathcal{R}(A+B).$$ By  \cite[Lemma~3.8]{LSX}, the inclusion  $\mathcal{R}(A)\subseteq \mathcal{R}(A+B)$ is valid whenever
$\mathcal{R}(A+B)$ is closed in $H$.
\end{remark}

\begin{theorem}\label{thm:unsolvable}There exist certain Hilbert space $H$ and $A,B\in\mathbb{B}(H)_+$ such that equation \eqref{equality is obtained} has no solution.
\end{theorem}
\begin{proof}Let $H=\ell^2(\mathbb{N})\oplus \ell^2(\mathbb{N})$ and let $A_0,B_0\in\mathbb{B}(H)$ be given by \eqref{counterexample A B}, where $S$ is defined by \eqref{operator s for A B}.
Clearly $B_0^\frac12=B_0$, and it is a routine matter to verify that
$$\left(
    \begin{array}{cc}
      1 & t \\
      t & t^2 \\
    \end{array}
  \right)^{\frac12}=\frac{1}{\sqrt{1+t^2}}\left(
    \begin{array}{cc}
      1 & t \\
      t & t^2 \\
    \end{array}
  \right),
  \quad t\in [0,+\infty).$$
Hence,
$$A_0^\frac12=\left(
                \begin{array}{cc}
                  (I+S^2)^{-\frac12} &\quad (I+S^2)^{-\frac12}S \\
                 (I+S^2)^{-\frac12}S &\quad  (I+S^2)^{-\frac12}S^2 \\
                \end{array}
              \right).$$
So, we may use the injectivity of $S$ to conclude that
$$\mathcal{R}\big(A_0^\frac12\big)\cap \mathcal{R}\big(B_0^\frac12\big)=(0,0)^T\in H.$$
Since $A_0:B_0\ge 0$ and  $\mathcal{R}\big[(A_0:B_0)^\frac12\big]=\mathcal{R}\big(A_0^\frac12\big)\cap \mathcal{R}\big(B_0^\frac12\big)$ \cite[Theorem~4.2]{Fillmore-Williams}, we see that $A_0:B_0=0$. Therefore, if we take $A=A_0$ and $B=B_0$, then equation
 \eqref{equality is obtained} is simplified as
\begin{equation*}\label{equality is obtained-special}
  0=X^*A_0X+(I-X)^*B_0(I-X), \quad X\in\mathbb{B}(H).
\end{equation*}
Suppose on the contrary that the above equation has a solution $X\in\mathbb{B}(H)$.
Then by the positivity of $A_0$ and $B_0$, we have
\begin{equation}\label{two zeros-1}B_0(I-X)=A_0X=0.\end{equation}
Let $X$ be partitioned as
$$X=\left(
                       \begin{array}{cc}
                         X_1 & X_2 \\
                         X_3 & X_4 \\
                       \end{array}
    \right),$$
in which $X_i\in\mathbb{B}\big(\ell^2(\mathbb{N})\big)$.  From \eqref{counterexample A B} and \eqref{two zeros-1}, it is easily seen that
$$X_1=I,\quad X_2=0,\quad I+SX_3=0,$$
which gives $S(-X_3)=I$. So, the same contradiction also occurs,  as shown at the end of the proof of Theorem~\ref{thm:solvability does not imp reduced solu}.
\end{proof}

We end this paper with the following concluding remark.
\begin{remark}Four equations are involved  in Definition~\ref{defn:(S,T)-weakly} for an $(M,N)$-weakly complementable operator.
We are unknown whether these equations can be simplified into the second and the third equations, as in the Hilbert space case.
\end{remark}

%


\vspace{5ex}

\begin{thebibliography}{99}


\bibitem{Anderson}W. N. Anderson, Jr., Shorted operators, SIAM J. Appl. Math. 20 (1971), 520--525.

\bibitem{Anderson-Duffin} W. N. Anderson, Jr. and R. J. Duffin, Series and parallel addition of matrices, J. Math. Anal. Appl. 26 (1969), 576--594.



\bibitem{AT}W. N. Anderson, Jr. and G. E. Trapp,  Shorted operators. II, SIAM J. Appl. Math. 28 (1975), 60--71.



\bibitem{ACS}J. Antezana, G. Corach and  D. Stojanoff, Bilateral shorted operators and parallel sums, Linear Algebra Appl. 414 (2006), no. 2--3, 570--588.


\bibitem{Douglas}R. G. Douglas, On majorization, factorization, and range inclusion of operators on Hilbert
space, Proc. Amer. Math. Soc. 17 (1966), 413--415.

\bibitem{Fang-Moslehian-Xu}X. Fang, M. S. Moslehian and Q. Xu, On majorization and range inclusion of operators on Hilbert $C^*$-modules, Linear Multilinear Algebra 66 (2018), no. 12, 2493--2500.



\bibitem{Fillmore-Williams}P. A. Fillmore and J. P. Williams, On operator ranges, Adv. Math. 7 (1971), 254--281.

\bibitem{Frank} M. Frank, Regularity results for classes of Hilbert $C^*$-modules with respect to special bounded modular functionals, Ann. Funct. Anal. 15 (2024), no. 2,  Paper No. 19, 18 pp.

\bibitem{Frank02}M. Frank, Hilbert $C^*$-modules over monotone complete $C^*$-algebras, Math. Nachr. 175 (1995), 61--83.  



\bibitem{FMXZ}C. Fu, M. S. Moslehian, Q. Xu and A. Zamani, Generalized parallel sum of adjointable operators on Hilbert $C^*$-modules, Linear Multilinear Algebra 70 (2022), no. 12, 2278--2296.


\bibitem{Hansen}F. Hansen, A note on the parallel sum, Linear Algebra Appl. 636 (2022), 69--76.




\bibitem{KS} J. Kaad and M. Skeide, Kernels of Hilbert module maps: a counterexample, J. Operator Theory 89 (2023), no. 2, 343--348. 




\bibitem{Lance}E. C. Lance, Hilbert $C^*$-modules--A toolkit for operator algebraists, Cambridge University Press, Cambridge, 1995.




\bibitem{LLX}N. Liu, W. Luo and Q. Xu, The polar decomposition for adjointable operators on Hilbert $C^*$-modules and centered operators, Adv. Oper. Theory 3 (2018), no. 4, 855--867.



\bibitem{LSX}W. Luo, C. Song and Q. Xu, The parallel sum for adjointable operators on Hilbert $C^*$-modules, (Chinese) Acta Math. Sinica (Chinese Ser.) 62 (2019), no. 4, 541--552.


\bibitem{Manuilov}V. M. Manuilov, On extendability of functionals on Hilbert $C^*$-modules, Math. Nachr. 297 (2024), 998--1005.

\bibitem{MMX}V. Manuilov, M. S.  Moslehian and Q. Xu, Douglas factorization theorem revisited, Proc. Amer. Math. Soc. 148 (2020), no. 3, 1139--1151.

\bibitem{MT} V. M. Manuilov and E. V. Troitsky, Hilbert $C^*$-modules, Translated from the 2001 Russian original by the authors, Translations of Mathematical Monographs, 226, American Mathematical Society, Providence, RI, 2005.

\bibitem{MT02} V. M. Manuilov and E. V. Troitsky, Hilbert $C^*$-modules with Hilbert dual and $C^*$-Fredholm operators, Integral Equations Operator Theory 95 (2023), no. 3, Paper No. 17, 12 pp.

\bibitem{Mitra-Odell}S. K. Mitra and P. L. Odell, On parallel summabillty of matrices, Linear  Algebra Appl. 74 (1986), 239--255.

\bibitem{Morley}T. D. Morley, An alternative approach to the parallel sum, Adv. Appl. Math. 10 (1989), 358--369.




\bibitem{Paschke}W. L. Paschke, Inner product modules over $B^*$-algebras, Trans. Amer. Math. Soc. 182 (1973), 443--468.


\bibitem{Pedersen} G. K. Pedersen, $C^*$-algebras and their automorphism groups, Academic Press, New York, 1979.

\bibitem{SMXZ}M. Sababheh, H. R. Moradi, Q. Xu and S. Zhao, Some inequalities for adjointable operators in Hilbert $C^*$-modules, arXiv:2402.13479

\bibitem{TWD}X. Tian, S. Wang and C. Deng, On parallel sum of operators, Linear Algebra Appl. 603 (2020), 57--83.

\bibitem{VMX}M. Vosough, M. S.  Moslehian and Q. Xu, Closed range and nonclosed range adjointable operators on Hilbert $C^*$-modules, Positivity 22 (2018), no. 3, 701--710.


\bibitem{Xu-Wei-Gu}Q. Xu, Y. Wei and Y. Gu, Sharp norm estimations for Moore-Penrose inverses of stable
perturbations of Hilbert $C^*$-module operators, SIAM J. Numer. Anal. 47 (2010), 4735--4758.






\end{thebibliography}
\end{document}